\documentclass{amsart}
\usepackage[dvips]{epsfig}
\usepackage{graphicx}
\usepackage{latexsym}
\usepackage{amsmath}
\usepackage{amsthm}
\usepackage{amssymb}
\usepackage{mathrsfs}
\usepackage[shortlabels]{enumitem}
\newlist{propenum}{enumerate}{1} 
\setlist[propenum]{label=(\roman*)}

\RequirePackage[colorlinks=true,citecolor=blue,urlcolor=blue]{hyperref}

\usepackage{caption}
\usepackage{subcaption}

\usepackage{eucal}

\usepackage{xcolor,soul}

\usepackage{mathtools}
\newcommand{\widesim}[3][1.5]{
  \mathrel{\underset{#3}{\overset{#2}{\scalebox{#1}[1]{$\sim$}}}}
}

 \usepackage[applemac]{inputenc}

\newtheorem{thm}{Theorem}[section]

\newtheorem{lem}[thm]{Lemma}
\newtheorem{cor}[thm]{Corollary}
\newtheorem{defi}[thm]{Definition}
\newtheorem{hyp}[thm]{Assumption}

\newtheorem{prop}[thm]{Proposition}
\theoremstyle{remark}

\newtheorem{rem}[thm]{Remark}

\newcommand{\zz}{\boldsymbol{z}}

\newcommand{\vip}{\vskip.2cm}
\newcommand{\COMMENTAIRE}[1]{}

\newcommand{\field}[1]{\mathbb{#1}}
\newcommand{\II}{\field{I}}
\newcommand{\JJ}{\field{J}}
\newcommand{\EE}{\field{E}}

\newcommand{\GG}{\field{G}}

\newcommand{\NN}{\field{N}}
\newcommand{\PP}{\field{P}}

\newcommand{\RR}{\field{R}}
\newcommand{\TT}{\field{T}}

\newcommand{\XX}{\field{X}}

\newcommand{\Bb}{{\mathcal B}}
\newcommand{\Cc}{{\mathcal C}}

\newcommand{\Ff}{{\mathcal F}}

\newcommand{\Hh}{{\mathcal H}}

\newcommand{\Pp}{{\mathcal P}}

\newcommand{\Qq}{{\mathcal Q}}

\newcommand{\rd}{{\rm d}}

\newcommand{\bF}{{\mathfrak f}}

\newcommand{\cb}{{\mathcal B}}

\newcommand{\cf}{{\mathcal F}}

\newcommand{\cp}{{\mathcal P}}
\newcommand{\cq}{{\mathcal Q}}

\newcommand{\cs}{{\mathscr S}}

\newcommand{\A}{{\mathbb A}}

\newcommand{\E}{{\mathbb E}}

\newcommand{\G}{\mathbb{G}}
\newcommand{\N}{{\mathbb N}}

\newcommand{\R}{{\mathbb R}}

\newcommand{\T}{\mathbb{T}}

\newcommand{\ind}{{\bf 1}}

\newcommand{\sot}{\otimes_{\rm sym}}

\newcommand{\norm}[1]{\mathop{\parallel\! #1 \! \parallel}\nolimits}

\newcommand{\reff}[1]{(\ref{#1})}

\newcommand{\superexp} {\xRightarrow[b_{n}^{2}]{\rm superexp}}
\newcommand{\equiexp} {\widesim[2]{{\rm superexp}}{b_{n}^{2}}}

\begin{document}

\title[mdp for kernel density estimator of BMC]{Moderate deviation principles for kernel estimator of invariant density in bifurcating Markov chains models.}

\author{S. Val\`ere Bitseki Penda}

\address{S. Val\`ere Bitseki Penda, IMB, CNRS-UMR 5584, Universit\'e Bourgogne Franche-Comt\'e, 9 avenue Alain Savary, 21078 Dijon Cedex, France.}

\email{simeon-valere.bitseki-penda@u-bourgogne.fr}



\begin{abstract}
Bitseki and Delmas (2021) have studied recently the central limit theorem for kernel estimator of invariant density in bifurcating Markov chains models. We complete their work by proving a moderate deviation principle for this estimator. Unlike the work of Bitseki and Gorgui (2021), it is interesting to see that the distinction of the two regimes disappears and that we are able to get moderate deviation principle for large values of the ergodic rate. It is also interesting and surprising to see that for moderate deviation principle, the ergodic rate begins to have an impact on the choice of the bandwidth for values smaller than in the context of central limit theorem studied by Bitseki and Delmas (2021).  
\end{abstract}

\maketitle






\section{Introduction}
The study of bifurcating Markov chains (BMCs, for short) models has taken a special place in the literature these last years due to their links with the study of  the cell dynamics (see for e.g. \cite{BO17, DM10, dhkr:segf, Guyon, MR4025704}). The first model of BMC, named ``symmetric'' bifurcating autoregressive process (BAR, for short) were introduced by Cowan and Staudte \cite{CS86} in order to understand the cell division mechanisms of Escherichia Coli (E. Coli, for short). E. Coli is a rod shaped bacterium which reproduces by dividing in two, thus producing two new cells. One of type 1 which has the old end of the mother and the other of type 0 which has the new end of the mother. The age of a cell is thus given by the age of its old pole in the sense of the number of divisions from which this pole exists. This cell division mechanism raises several questions, among other that of the symmetry of the division. In order to give a rigorous answer to this question, Guyon \cite{Guyon} has developed and studied the theory of BMCs. We note that to the best of our knowledge, the term BMC appears for the first time in the work of \cite{BZ04}. In particular, Guyon has studied an extension of the model introduced by Cowan and Staudte, named ``asymmetric'' BAR. In the conclusion of his study, Guyon concludes that aging has an impact on cell reproduction. We note that an extension of the model proposed by Guyon, named nonlinear BAR (NBAR, for short) were studied by Bitseki and Olivier in \cite{BO17}. Another question of interest related to cell division is estimating the division rate at which cells divide. This question has been tackled recently in the work of Doumic \& \textit{al.} \cite{dhkr:segf} and Hoffman \& Marguet \cite{MR4025704}. In all the previous work, the behaviour and the definition of parameters of interest are associated with the density of the invariant probability of an auxiliary Markov chain (see below for a precise definition). The estimation of this invariant density has recently been the subject of several studies. One can cite \cite{BHO, bitsekiroche2020} where adaptive methods have been proposed for the estimation of this invariant density. More recently, Bitseki and Delmas \cite{BD3} have studied central limit theorem for kernel estimators of this invariant density. Our main objective in this paper is to complete the previous study by establishing a moderate deviation principle for these kernel estimators. Before going any further, let us recall the definition of the main concepts that we will use and study.

\section{The model of bifurcating Markov chain and definition of the estimators}\label{sec:model}

\subsection{The regular binary tree associated to BMC models}
We denote by $\NN$ (resp. $\NN^{*}$) the space of (resp. positive) natural integers. We   set   $\T_0=\G_0=\{\emptyset\}$, $\G_k=\{0,1\}^k$  and $\T_k  =  \bigcup _{0  \leq r  \leq  k} \G_r$  for $k\in  \N^*$, and  $\T  =  \bigcup _{r\in  \N}  \G_r$. The set  $\G_k$ corresponds to the  $k$-th generation, $\T_k$ to the tree  up to the $k$-th generation, and $\T$ the complete binary  tree. One can see that the genealogy of the cells is entirely described by $\TT$ (each vertex of the tree designates an individual). For $i\in \T$, we denote by $|i|$ the generation of $i$ ($|i|=k$  if and only if $i\in \G_k$) and $iA=\{ij; j\in A\}$  for $A\subset \T$, where $ij$  is the concatenation of   the  two   sequences  $i,j\in   \T$,  with   the  convention   that $\emptyset i=i\emptyset=i$. For $A \subset \TT$, we denote by $|A|$ the number of elements of $A$. Note that for all $n \in \NN,$ $|\GG_{n}| = 2^{n}$ and $|\TT_{n}| = 2^{n+1} - 1.$

\subsection{The probability kernels associated to BMC models} 
\quad\\
For our convenience, we set $S = \RR^{d}$, $d \geq 1$ and $S$ is equipped with the Borel sigma-algebra $\cs$. For any $q \in \NN^{*}$, we denote by $\Bb(S^{q})$ (resp. $\Bb_{b}(S^{q})$, resp. $\Cc_{b}(S^{q})$) the space of  (resp. bounded, resp. bounded continuous ) $\RR\text{-}$valued measurable functions defined on $S^{q}$. For all $q \in \NN^{*}$, we set $\cs^{\otimes q} = \cs \otimes \ldots \otimes \cs$.  Let $\Pp$ be a probability kernel on $(S,\cs^{\otimes2})$, that is: $\Pp(\cdot  , A)$  is measurable  for all  $A\in \cs^{\otimes 2}$,  and $\Pp(x, \cdot)$ is  a probability measure on $(S^2,\cs^{\otimes 2})$ for all $x \in S$. For any $g\in \cb_b(S^3)$ and $h\in \cb_b(S^2)$,   we set for $x\in S$:
\begin{equation}\label{eq:Pg}
(\Pp g)(x)=\int_{S^2} g(x,y,z)\; \Pp(x,\rd y,\rd z) \quad \text{and} \quad (\Pp h)(x)=\int_{S^2} h(y,z)\; \Pp(x,\rd y,\rd z).
\end{equation}
We define $(\Pp g)$ (resp. $(\Pp h)$), or simply $Pg$ for $g\in \cb(S^3)$(resp. $\Pp h$ for $h\in \cb(S^2)$), as soon as the corresponding integral \reff{eq:Pg} is well defined, and we have  that $\Pp g$ and $\Pp h$ belong to $\cb(S)$. We denote by $\Pp_{0}$, $\Pp_{1}$ and $\Qq$ respectively the first and the second marginal of $\Pp$, and the mean of $\Pp_{0}$ and $\Pp_{1}$, that is, for all $x \in S$ and $B \in \mathcal{S}$  
\begin{equation*}\label{eq:P0P1Q}
\Pp_{0}(x,B) = \Pp(x,B\times S), \quad \Pp_{1}(x,B) = \Pp(x,S\times B) \quad \text{ and} \quad \Qq = \frac{(\Pp_{0} + \Pp_{1})}{2}.
\end{equation*}

\medskip 

Now let us give  a precise definition of bifurcating Markov chain. 
\begin{defi}[Bifurcating Markov Chains, see \cite{Guyon, BD3}]
\quad\\
We say  a stochastic process indexed  by $\T$, $X=(X_i,  i\in \T)$, is a bifurcating Markov chain (BMC) on a measurable space $(S, \cs)$ with initial probability distribution  $\nu$ on $(S, \cs)$ and probability kernel $\cp$ on $S\times \cs^{\otimes 2}$ if:
\begin{itemize}
\item[-] (Initial  distribution.) The  random variable  $X_\emptyset$ is distributed as $\nu$.
\item[-] (Branching Markov property.) For  any  sequence   $(g_i, i\in \T)$ of functions belonging to $\cb_b(S^3)$ and  for all $k\geq 0$, we have
\begin{equation*}
\E\Big[\prod_{i\in \G_k} g_i(X_i,X_{i0},X_{i1}) |\sigma(X_j; j\in \T_k)\Big] 
=\prod_{i\in \G_k} \cp g_i(X_{i}).
\end{equation*}
\end{itemize}
\end{defi}
Following \cite{Guyon}, we introduce an  auxiliary Markov  chain $Y=(Y_n, n\in  \N) $  on $(S,\cs)$ with $Y_{0}=X_{1}$ and transition probability $\Qq$. The chain $(Y_{n}, n\in \mathbb{N})$ corresponds to a random lineage taken in the population. We  shall   write  $\E_x$   when  $X_\emptyset=x$ (\textit{i.e.}  the initial  distribution  $\nu$ is  the  Dirac mass  at $x\in S$). We will assume that the Markov chain $Y$ is ergodic and we denote by $\mu$ its invariant probability measure. Asymptotic and non-asymptotic behaviour of BMCs are strongly related to the knowledge of $\mu$. In particular, Guyon has proved that if $Y$ is ergodic, then for all $f \in \Cc_{b}(S)$,
\begin{equation*}
|\A_{n}|^{-1} \sum_{u \in \A_{n}} f(X_{u}) \underset{n \rightarrow +\infty}{\xrightarrow{\hspace*{0.75cm}}} \langle \mu,f \rangle \quad \text{in probability,} \quad \text{where $\A_{n} \in \{\GG_{n}, \TT_{n}\}.$}
\end{equation*}
But in most cases, the invariant probability $\mu$ is unknown, so its estimation from the data is of great interest. For that purpose, we do the following assumption.

\begin{hyp}\label{hyp:DenMu}
The transition kernel $\Pp$ has a density, still denoted by $\Pp$, with respect to the Lebesgue measure. 
\end{hyp}
\begin{rem}\label{rem:}
Assumption \ref{hyp:DenMu} implies that the transition kernel $\Qq$ has a density, still denoted by $\Qq$, with respect to  the Lebesgue measure. More precisely, we have $\Qq(x,y) = 2^{-1} \int_{S} (\Pp(x,y,z)+\Pp(x,z,y))dz.$ This implies in particular that the invariant probability $\mu$ has a density, still denoted by $\mu$, with respect to the Lebesgue measure (for more details, we refer for e.g. to \cite{duflo2013random}, chap 6).
\end{rem}

\subsection{Kernel estimator of the invariant density $\mu$}

Recall that $\A_{n} \in \{\GG_{n}, \TT_{n}\}$ and $S = \RR^{d}$, $d \geq 1$. Assume we observe $\XX^{n} = (X_{u}, u \in \A_{n})$. Let $(h_{n}, n \in \NN)$ be a sequence of positive numbers which converges to $0$ as $n$ goes to infinity. We will simply write $h$ for $h_{n}$ if there is no ambiguity.  Let the kernel function $K: S \rightarrow \RR$ such that $\int_{S} K(x)dx = 1.$ Then, for all $x \in S,$ we propose to estimate $\mu(x)$ by
\begin{equation}\label{eq:def-est}
\widehat{\mu}_{\A_{n}}(x) = |\A_{n}|^{-1} h_{n}^{-d/2} \sum_{u \in \A_{n}} K_{h_{n}}(x - X_{u}),
\end{equation}
where $K_{h_{n}}(\cdot) = h_{n}^{-d/2}K(h_{n}\cdot).$ These estimators are strongly inspired from \cite{masry1986, Parzen1962, Roussas1969a}.  They have been studied in \cite{dhkr:segf, bitsekiroche2020} (non asymptotic studies) and in \cite{BD3} (central limit theorem). 

\subsection{Moderate deviation principle and related topics}

Our aim is to study moderate deviation principles for the estimators defined in \eqref{eq:def-est}. Before we proceed, let us introduce the notion of moderate deviation principle. We give the definition in a general setting. Let $(Z_{n})_{n\geq 0}$ be a sequence of random variables with values in $S$ endowed with its Borel $\sigma$-field $\cs$ and let $(s_{n})_{n\geq 0}$ be a positive sequence that converges to $+\infty$. We assume that $Z_{n}/s_{n}$ converges in probability to 0 and that $Z_{n}/\sqrt{s_{n}}$ converges in distribution to a centered Gaussian law. Let $I: S \rightarrow\RR^{+}$ be a lower semicontinuous function, that is for all $c>0$ the sub-level set $\{x \in S, I(x)\leq c\}$ is a closed set. Such a function $I$ is called {\it rate function} and it is called {\it good rate function} if all its sub-level sets are compact sets. Let $(b_{n})_{n\geq 0}$ be a positive sequence such that $b_n \rightarrow + \infty $ and $b_{n}/\sqrt{s_{n}} \rightarrow 0$ as $n$ goes to $+\infty$.
\begin{defi}[Moderate deviation principle, MDP] \label{def:mdp}
	\quad\\
We say that $Z_{n}/(b_{n}\sqrt{s_{n}})$ satisfies a moderate deviation principle on $S$ with speed $b_{n}^{2}$ and rate function~$I$ if, for any $A\in\cs$,
\begin{align*}
-\inf_{x\in \mathring{A}}I(x) \leq \liminf_{n\rightarrow+\infty}\frac{1}{b_{n}^{2}} \log\PP\big(\frac{Z_{n}}{b_{n}\sqrt{s_{n}}}\in A\big)  \leq \limsup_{n\rightarrow+\infty}\frac{1}{b_{n}^{2}} \log\PP\big(\frac{Z_{n}}{b_{n}\sqrt{s_{n}}}\in A\big) \leq -\inf_{x\in\bar{A}}I(x),
\end{align*}
where $\mathring{A}$ and $\bar{A}$ denote respectively the interior and the closure of $A$.
\end{defi}

The following two concepts are closely related to the theory of  MDP: super-exponential convergence and exponential equivalence. Let $(Z_{n}, n \in \NN)$, $(W_{n}, n \in \NN)$ be sequences of random variables and $Z$ a random variable with value in a metric space $(S,d)$. 
\begin{defi}[Super-exponential convergence] \label{def:cv-super}
We say that $(Z_{n})_{n \geq 0}$ converges \, $(b_{n}^2)\text{-}$ super-exponentially fast in probability to $Z$ and we note $Z_{n} \xRightarrow[b_{n}^{2}]{\rm superexp} Z$ if, for all $\delta > 0$,
\begin{equation*}
\limsup_{n \rightarrow +\infty} \frac{1}{b_{n}^2} \log \PP\big( d( Z_{n}, Z) > \delta \big) = -\infty.
\end{equation*}
\end{defi}

\begin{defi}[Exponential equivalence, see \cite{DZ1998}, Chap 4] \label{def:equi-expo}
We say that $(Z_{n})_{n\geq 0}$ and $(W_{n})_{n\geq 0}$ are $(b_{n}^{2})_{n\geq 0}$-exponential\-ly equivalent and we note $Z_{n} \widesim[2]{{\rm superexp}}{b_{n}^{2}}  W_{n}$ if for any $\delta>0$,
$$
\limsup_{n\rightarrow+\infty}\frac{1}{b_{n}^{2}}
\log\PP\big(d(Z_{n},W_{n}) > \delta\big) = -\infty.
$$
\end{defi} 
 
\begin{rem}\label{rem:cv-det-exp} 
Note that for a determininistic sequence that converges to some limit $\ell$, it also converges $(b_{n}^2)\text{-}$superexponentially fast to $\ell$ for any rate $b_n$. We also note that if $(Z_{n})_{n\geq 0}$ and $(W_{n})_{n\geq 0}$ are $(b_{n}^{2})_{n\geq 0}$-exponential\-ly equivalent and if $(Z_{n})_{n\geq 0}$ satisfies a MDP, then $(W_{n})_{n\geq 0}$ satisfies the same MDP (for more details, see for e.g \cite{DZ1998}, Chap 4).
\end{rem}
The following result give a sufficient condition for super-exponential convergence of a sequence of random variables.	
\begin{rem}\label{rem3}
We assume that $(S,d)$ is a metric space. Let  $\left(Z_n\right)_{n \in \NN}$ be a sequence of random variables with values in $S$,  $Z$ a random variable with values in $S$. So if  $d(Z_n,Z) $ is  upper-bounded by  a deterministic  sequence which converges to  $0$,  then, for all sequence $(b_{n}, n \in \NN)$ converging to $+\infty$, $Z_{n} \xRightarrow[b_{n}^{2}]{\rm superexp} Z$.
\end{rem}

%
%

The moderate deviation principle has been proved in the i.i.d. setting for kernel density estimator, see for e.g. Gao \cite{gao2003moderate}, Mokkadem \& \textit{al.} \cite{mokkadem2005large}. We refer also to \cite{mokkadem2006confidence} where Mokkaddem and Pelletier have constructed confidence bands for probability densities based on moderate deviation principles.  In this paper, we will establish moderate deviation principle for $\widehat{\mu}_{\A_{n}}(x)$ following the martingale approach developed in \cite{BD3}. We will need the following assumption.
\begin{hyp}\label{hyp:erg-unif}
There exists a positive real number $M$  and $\alpha\in (0, 1)$ such that for all $f\in  \Bb_{b}(S)$: 
\begin{equation}\label{eq:geom-erg}
|\cq^{n}f - \langle \mu, f \rangle| \leq M \,  \alpha^{n} \|f\|_{\infty} \quad \text{for all  $n\in \N$.}
\end{equation}
\end{hyp}
\begin{rem}
Assumption \ref{hyp:erg-unif} is for example satisfy for nonlinear bifurcating autoregressive process under mild hypotheses on the autoregression functions (see \cite{BO2018} Lemma 9 for more details).
\end{rem}
The others assumptions we will need are based on the following bias-variance type decomposition of the estimator $\widehat{\mu}_{\A_{n}}(x)$: 
\begin{equation}\label{eq:DeBiVa}
\widehat{\mu}_{\A_{n}}(x) - \mu(x) = B_{h_n}(x) + V_{\A_n, h_n}(x),
\end{equation}
where for $h>0$ and $\A\subset \T$ finite:
\[
B_h(x)=h^{-d/2}K_{h}\star\mu(x) - \mu(x) \quad\text{and}\quad V_{\A, h}(x)=  |\A|^{-1}h^{-d/2} \sum_{u \in \A} \Big(K_{h}(x - X_{u}) - K_{h}\star\mu(x)\Big),
\]
and for $h>0$ and $u\in \T$, we set:
\[
  K_{h}\star\mu(x) = \EE_{\mu}[K_{h}(x-X_{u})] = \int_{S}
K_{h}(x-y) \mu(y)\, dy.
\]
To study the variance term $V_{\A_{n},h_{n}}(x)$, we will introduce a more general sequence of functions (see Section \ref{sec:mdp-ad}). 

The following assumptions on the kernel, the bandwidth and the regularity of the unknown density function are usual. Recall  $S=\R^d$ with $d\geq 1$.
\begin{hyp}[Regularity of the kernel function and the bandwidth]\label{hyp:K}$\,$
\begin{propenum}
\item\label{item:cond-f}  The kernel function $K\in \cb(S)$  satisfies:
\begin{equation*}
\norm{K}_{\infty}<+\infty, \,\, \norm{K}_1      < + \infty, \,\, \norm{K}_{2} <+\infty, \,\, \int_{S}\!  K(x)\,   dx  =  1 \quad\text{and}\quad \lim_{|x|\rightarrow  +\infty}       |x|K(x)=0.
\end{equation*} 
\item\label{item:bandwidth} There exists $\gamma\in (0, 1/d)$ such that the    bandwidth  $(h_n,n \in \NN)$  are  defined by $h_n= 2^{-n\gamma}$. 
\end{propenum}
\end{hyp}

\begin{hyp}[Further regularity on the density $\mu$, the kernel function and the bandwidths]\label{hyp:estim-tcl}

Suppose that there exists an invariant probability measure $\mu$  of $\cq$ and that Assumptions \ref{hyp:DenMu} and \ref{hyp:K} hold. We assume there exists  $s > 0$ such that the following hold:
\begin{propenum}
\item\label{item:den-Holder}  \textbf{The density  $\mu$ belongs  to the (isotropic) H\"older class  of order $(s, \ldots,  s) \in \RR^{d}$:} The  density $\mu$ admits partial derivatives with respect to  $x_{j}$, for all  $j\in \{1,\ldots  d\}$, up to  the order $\lfloor s  \rfloor$ and there exists  a finite constant $L  > 0$ such that  for all  $x=(x_1,  \ldots,  x_d), \in  \RR^{d}$,  $t\in \R$  and $ j \in \{1, \ldots, d\}$:
\begin{equation*}
\left|\frac{\partial^{\lfloor s \rfloor}\mu}{\partial x_{j}^{\lfloor s \rfloor}}(x_{-j},t)-\frac{\partial^{\lfloor s \rfloor}\mu}{\partial x_{j}^{\lfloor s \rfloor}}(x)\right| \leq L|x_{j} - t|^{\{s\}},  
\end{equation*}
where $(x_{-j},t)$ denotes the vector $x$ where we have replaced the $j^{th}$ coordinate $x_{j}$ by $t$, with the convention ${\partial^{0}\mu}/{\partial x_{j}^{0}} = \mu$. 

\item\label{item:K}\textbf{The kernel $K$ is  of order $(\lfloor s \rfloor, \ldots, \lfloor s \rfloor) \in \NN^{d}$:} We have $\int_{\RR^{d}} |x|^{s}K(x)\, dx < \infty$  and $\int_{\RR} x^{k}_{j}\, K(x)\, dx_{j} = 0$ for all $k \in \{1,\ldots,\lfloor s \rfloor\}$ and $j \in \{1,\ldots,d\}$. 

\end{propenum}
\end{hyp}
For $\alpha > 1/2$, we shall also assume the following.
\begin{hyp}\label{ass:fln-aS2}
Keeping the same notations as in \ref{item:bandwidth} of Assumption \ref{hyp:K}, we further assume that Assumption \ref{hyp:erg-unif} holds with
\begin{equation}\label{eq:fln-aS2}
\lim_{n \rightarrow +\infty} (2^{1-d\gamma}\alpha)^{n} = 0.
\end{equation}
\end{hyp}
\begin{rem}
As consequence of Assumption \ref{ass:fln-aS2} and \ref{item:bandwidth} of Assumption \ref{hyp:K}, for moderate deviation principle, the ergodicity rate $\alpha$ begins to have an impact on the choice of the bandwidth for $\alpha > 1/2.$ This is out of step with the central limit theorem where the ergodicity rate $\alpha$ begins to have an impact on the choice of the bandwidth for $\alpha > 1/\sqrt{2}$ (see \cite{BD3} for more details).
\end{rem}

In the sequel, we will consider the positive sequence
$(b_{n},n \in \NN)$ such that:
\begin{equation}\label{eq:speed-mdp}
\lim_{n \rightarrow +\infty} b_{n} = +\infty; \quad \lim_{n \rightarrow +\infty} \frac{n^{3/2} \, b_{n}}{\sqrt{|\GG_{n}|h_{n}^{d}}} = 0; \quad \lim_{n \rightarrow +\infty} \frac{b_{n}}{\sqrt{|\GG_{n}| h_{n}^{2s + d}}} = + \infty,
\end{equation}
where $s$ is the regularity parameter given in Assumption \ref{hyp:estim-tcl}.

The paper is organised as follows. In Section \ref{sec:MDP-m-SEC} we state the main result for the moderate deviation principles of the estimators $\widehat{\mu}_{\A_{n}}(x)$ for $x$ in the set continuity of $\mu$ and $\A_{n} \in \{\TT_{n}, \GG_{n}\}$. In Section \ref{sec:mdp-ad}, directly linked to the study of variance term $V_{\A,h}(x)$ defined in \eqref{eq:DeBiVa}, we study the moderate deviation principle for general additive functionals of BMCs. Sections \ref{sec:proof-main1} and \ref{sec:proof-main2} are devoted to the proofs of results. In Section \ref{sec:appendix}, we recall some useful results.

\section{Main result}

\subsection{Moderate deviation principle for $\widehat{\mu}_{\A_{n}}$}\label{sec:MDP-m-SEC}

First, we state a strong consistency result for the estimators $\widehat{\mu}_{\A_{n}}(x)$ for $x$ in the set of continuity of $\mu$. Its proof is given in Section \ref{sec:Pstr-cv-mu}.
\begin{lem}\label{lem:str-cv-mu}
Let $X$ be a BMC with kernel $\Pp$ and initial distribution $\nu$ such that Assumptions \ref{hyp:erg-unif}, \ref{hyp:K} and \ref{hyp:estim-tcl} hold. Furthermore, if $\alpha > 1/2$ then assume that Assumption \ref{ass:fln-aS2} holds. Let $(b_{n}, n \in \NN)$ be a positive sequence with satisfies \eqref{eq:speed-mdp}. Then, for all $x$ in the set of continuity of $\mu$ and $\A_{n} \in \{\TT_{n},\GG_{n}\}$ we have $\widehat{\mu}_{\A_{n}}(x) \xRightarrow[b_{n}^{2}]{\rm superexp} \mu(x).$ 
\end{lem}

The main result of this Section is the following theorem which state the moderate deviation principle for $\widehat{\mu}_{\A_{n}}(x) - \mu(x)$ for $x$ in the set of continuity of the function $\mu$.
\begin{thm}\label{thm:m-est-main}
Under the hypothesis of Lemma \ref{lem:str-cv-mu}, for all $x$ in the set of continuity of $\mu$ and $\A_{n} \in \{\TT_{n}, \GG_{n}\}$, $b_{n}^{-1} \sqrt{|\A_{n}|h_{n}^{d}} (\widehat{\mu}_{\A_{n}}(x) - \mu(x))$ satisfies a moderate deviation principle on $\RR$ with speed $b_{n}^{2}$ and rate function $I$ defined by: $I(y) = y^{2}/(2 \|K\|_{2}^{2} \mu(x))$ for all $y \in \RR$, 
that is, for any $A \subset \RR$,
\begin{align*}
-\inf_{y\in \mathring{A}}I(y) \leq \liminf_{n\rightarrow+\infty}\frac{1}{b_{n}^{2}} &\log\PP\big(b_{n}^{-1} \sqrt{|\A_{n}|h_{n}^{d}} (\widehat{\mu}_{\A_{n}}(x) - \mu(x)) \in  A \big) \\  
&\leq \limsup_{n\rightarrow+\infty}\frac{1}{b_{n}^{2}} \log\PP\big(b_{n}^{-1} \sqrt{|\A_{n}|h_{n}^{d}} (\widehat{\mu}_{\A_{n}}(x) - \mu(x)) \in  A \big) \leq -\inf_{y \in \bar{A}}I(y),
\end{align*}
where  $\mathring{A}$ and $\bar{A}$ denote respectively the interior and the closure of $A$.
\end{thm}

In order to obtain confidence intervals for $\mu(x)$, it would be interesting to replace $\mu(x)$ in the expression of the rate function $I(\cdot)$ by an estimator. In that direction, we have the following. Let $\A_{n}^{*} \in \{\GG_{n}, \TT_{n}\}.$ Obviously, $\A_{n}^{*}$ and $\A_{n}$ can be the same. We consider the estimator $\widehat{\mu}_{\A_{n}^{*}}(x)$ of $\mu(x)$ defined with $\A_{n}^{*}$ instead of $\A_{n}$. Let $(\varpi_{n}, n \in \NN)$ be a sequence of real numbers such that $\varpi_{n} \rightarrow 0$ as $n \rightarrow +\infty.$ Then, we have the following result which the proof is given in Section \ref{sec:Pest-main2}.

\begin{thm}\label{thm:est-main2}
Under the hypothesis of Lemma \ref{lem:str-cv-mu}, for all $x$ in the set of continuity of $\mu$ and $\A_{n}, \A_{n}^{*} \in \{\TT_{n}, \GG_{n}\}$, $b_{n}^{-1} (\|K\|_{2} \sqrt{\widehat{\mu}_{\A_{n}^{*}}(x)} \vee \varpi_{n})^{-1} \sqrt{|\A_{n}|h_{n}^{d}} (\widehat{\mu}_{\A_{n}}(x) - \mu(x))$ satisfies a moderate deviation principle on $\RR$ with speed $b_{n}^{2}$ and rate function $I'$ defined by: $I'(y) = y^{2}/2$ for all $y \in \RR$. 
\end{thm}

In particular, using the contraction principle (see for e.g Dembo and Zeitouni \cite{DZ1998}, Chap 4), we have the following corollary of Theorem \ref{thm:est-main2}.
\begin{cor}\label{cor:est-main2} 
Under the hypothesis of Theorem \ref{thm:est-main2}, we have the following convergence for $x$ in the set of continuity of $\mu$ and $\A_{n}, \A^{*} \in \{\TT_{n}, \GG_{n}\}:$
\[
\lim_{n \rightarrow +\infty} \frac{1}{b_{n}^{2}} \log \PP\Big(b_{n}^{-1} \, \Big(\|K\|_{2} \sqrt{\widehat{\mu}_{\A_{n}^{*}}(x)} \vee \varpi_{n}\Big)^{-1} \, \sqrt{|\A_{n}|h_{n}^{d}} \Big|\big(\widehat{\mu}_{\A_{n}}(x) - \mu(x)\big)\Big| > \delta \Big) = - \frac{\delta^{2}}{2} \quad \forall \delta > 0.
\]
\end{cor}
\begin{rem}
Corollary \ref{cor:est-main2} yields a simple confidence interval for $\mu(x)$, of decreasing size \, \, $b_{n}/\sqrt{|\A_{n}| \, h_{n}^{d}}$ and with level asymptotically close to $1 - \exp(- (b_{n}^{2} \, \delta^{2})/2).$ 
\end{rem}

Using the structure of the asymptotic variance $\sigma^{2}$ in \eqref{eq:limf-ln}, we can prove the following multidimensional result which the proof is given in Section \ref{sec:est-mult}    
\begin{cor}\label{cor:m-est-mult}
Under the hypothesis of Theorem \ref{thm:m-est-main}, we have, for $x$ in the set of continuity of $\mu$ and for all $k \geq 0$, $b_{n}^{-1} \Big( |\GG_{n}|^{1/2}h_{n}^{1/2}\big(\widehat{\mu}_{\GG_{n}}(x) - \mu(x)\big), \ldots, |\GG_{n-k}|^{1/2}h_{n-k}^{1/2}\big(\widehat{\mu}_{\GG_{n-k}}(x) - \mu(x)\big) \Big)^{t}$ satisfies a moderate deviation principle on $\RR^{k+1}$ with speed $b_{n}^{2}$ and good rate function $J_x: \RR^{k+1} \rightarrow \RR$ defined by
\begin{equation*}
J_x(\zz) = \big( 2 \, \| K \|_2^2 \, \mu(x)\big)^{-1} \zz^{t}\Gamma^{-1} \zz \, , \quad \zz \in \RR^{k+1},
\end{equation*}
with $\Gamma = diag(2^{0}, \ldots, 2^{k})$, where $diag(\cdot)$ denotes the diagonal matrix and $\zz^{t}$ stands for the transpose of vector $\zz$.
\end{cor}
\begin{rem}
We deduce from Corollary \ref{cor:m-est-mult} that the estimators $|\G_{n-\ell}|^{1/2} h_{n-\ell}^{d/2} (\widehat{\mu}_{\G_{n-\ell}}(x) - \mu(x))$ are asymptotically independent in the sense of moderate deviation for $\ell \in \{0, \ldots, k\}$ and for any $k \in \NN.$
\end{rem}

\subsection{Moderate deviation principle for additive functionals of BMCs}\label{sec:mdp-ad}$\,$

 In order to study the variance term $V_{\A_{n},h_{n}}(x)$, we give here a moderate deviation principle for a general additive functionals of BMCs. For that purpose, we introduce the following assumption.
\begin{hyp}\label{ass:fl-n}
For $n \in \NN$, let $\bF_{n} = (f_{\ell,n}, n \geq \ell \geq 0)$ be a sequence of functions defined on $S$ such that $f_{\ell,n} = 0$ if $\ell > n$ and there exists $\gamma \in (0,1/d)$ such that:
\begin{itemize}
\item[(i)] $\sup_{0 \leq \ell \leq n} \{2^{-d \gamma n/2} \|f_{\ell,n}\|_{\infty}; \, 2^{d \gamma n/2} \|\cq f_{\ell,n}\|_{\infty}; \, \|\cq (f_{\ell,n}^{2})\|_{\infty}; \,  2^{d \gamma n} \|\Pp(f_{\ell,n} \otimes^{2}) \|_{\infty}\} < +\infty.$
\item[(ii)] $\sup_{0 \leq \ell \leq n} \{2^{d \gamma n/2}\langle \mu, |f_{\ell,n}| \rangle; \, \langle\mu, f_{\ell,n}^{2} \rangle\} < +\infty.$
\item[(iii)] The following limit exists and is finite:
\begin{equation}\label{eq:limf-ln}
\sigma^2 = \lim_{n\rightarrow+\infty } \sum_{\ell=0}^n 2^{-\ell} \norm{f_{\ell,n} }_{L^2(\mu)} ^2 < +\infty.
\end{equation} 
\end{itemize}
\end{hyp}
We will use the following notations. For a finite set $\A\subset \T$ and a function $f\in \cb(S)$, we set:
\begin{equation*}
M_\A(f)=\sum_{i\in \A} f(X_i).
\end{equation*}
In this paper, we are interested in the cases $\A = \GG_{n}$ and $\A = \TT_{n}$, that is the $n\text{-}$th generation and the first $n$ generation of the tree. Recall $\mu$ the invariant probability of $\Qq$, transition probability of the auxiliary Markov chain $(Y_{n}, n \in \NN)$. For $f \in L^{1}(\mu)$, we set: 
\begin{equation*}
\tilde{f} = f - \langle \mu, f\rangle.
\end{equation*}
Recall the sequence $\bF_{n}$ defined in Assumption \ref{ass:fl-n}. For $n \in \NN$, we set:
\begin{equation}\label{eq:Nemptyf}
N_{n,\emptyset}(\bF_{n}) = |\GG_{n}|^{-1/2} \sum_{\ell = 0}^{n} M_{\GG_{n-\ell}}(\tilde{f}_{\ell,n}).
\end{equation} 
The notation $N_{n,\emptyset}$ means that we consider the average from the root $\emptyset$ to the $n\text{-}$th generation.
\begin{rem}\label{rem:cas-P-Nnfn}
The definition of $N_{n,\emptyset}(\bF_{n})$ in \eqref{eq:Nemptyf} is mainly motivated by the decomposition \eqref{eq:DeBiVa}. It will allow us to threat the variance term of the estimator $\widehat{\mu}_{\A_{n}}(x)$ defined in \eqref{eq:def-est}. Instead, for $n \in \NN$, we set $f_{n}^{x}(\cdot) = K_{h_{n}}(x - \cdot)$.
Then, we consider the sequences of functions $(f_{\ell,n}^{\text{id}},\,  n  \geq \ell  \geq  0)$ and $(f^0_{\ell,n},\, n\geq \ell\geq 0)$ defined by:
\begin{equation}\label{eq:def-f-kernel}
f_{\ell,n}^{\text{id}} = f^{x}_{n} \quad \text{and} \quad f^{0}_{\ell, n} = f^x_n \ind_{\{\ell=0\}} . 
\end{equation}
It is not difficult to check that under Assumption \ref{hyp:K}, the sequence $(f_{\ell,n}^{\text{id}},\,  n  \geq \ell  \geq  0)$ and $(f^0_{\ell,n},\, n\geq \ell\geq 0)$ defined in \eqref{eq:def-f-kernel} satisfy Assumption \ref{ass:fl-n}. In particular, let $x$ be in the set of continuity of $\mu$.  Thanks to Lemma \ref{lem:bochner}, we have:
\begin{equation}\label{eq:lim-fnx-sig}
\lim_{n \rightarrow +\infty}  \norm{ f^x_{n}}_{L^2(\mu)}^2 = \lim_{n \rightarrow +\infty }  \langle \mu, (f^x_{n})^2  \rangle = \mu(x)\norm{K}_2^2. 
\end{equation}

If $\A_{n} = \GG_{n}$, it suffices to consider the sequence $\bF_{n} = (f_{\ell,n}, 0 \leq \ell \leq n)$ with $f_{\ell,n} = f_{\ell,n}^{0}$ and in that case, using \eqref{eq:lim-fnx-sig}, the asymptotic variance defined in \eqref{eq:limf-ln} is given by $\sigma^{2} = \|K\|_{2}^{2} \, \mu(x)$. If $\A_{n} = \TT_{n}$, it suffices to consider the sequence $\bF_{n} = (f_{\ell,n}, 0 \leq \ell \leq n)$ with  $f_{\ell,n} = f_{\ell,n}^{id}$ and in that case, using \eqref{eq:lim-fnx-sig}, the asymptotic variance defined in \eqref{eq:limf-ln} is given by $\sigma^{2} = 2 \|K\|_{2}^{2} \, \mu(x)$. 
\end{rem}
For our convenience, we assume that the quantity $\gamma$ which appears in Assumptions \ref{hyp:K} and \ref{ass:fl-n} is the same. The main result of this section is the following.
\begin{thm}\label{thm:mdp-fln}
Let $X$ be a BMC with kernel $\Pp$ and initial distribution $\nu$ such that Assumptions \ref{hyp:erg-unif}, \ref{hyp:K} and \ref{ass:fl-n} hold. Furthermore, if $\alpha > 1/2$ then assume that Assumption \ref{ass:fln-aS2} holds. Let $(b_{n}, n \in \NN)$ be a positive sequence with satisfies \eqref{eq:speed-mdp}. Then $b_{n}^{-1} N_{n,\emptyset}(\bF_{n})$ satisfies a moderate deviation principle on $\RR$ with speed $b_{n}^{2}$ and rate function $I$ defined by: $I(x) = x^{2}/(2 \sigma^{2})$ for all $x \in \RR$, with the finite variance $\sigma^{2}$ defined in \eqref{eq:limf-ln}.  
\end{thm}
\begin{rem}
In particular, using the contraction principle (see for e.g Dembo and Zeitouni \cite{DZ1998}, Chap 4), Theorem \ref{thm:mdp-fln} implies that
\begin{equation*}
\lim_{n \rightarrow +\infty} \frac{1}{b_{n}^{2}} \log \PP\left(\left|b_{n}^{-1} N_{n,\emptyset}(\bF_{n})\right| > \delta \right) = -I(\delta) \quad \forall \delta > 0.
\end{equation*}
\end{rem}
\begin{rem}
Unlike the results of Bitseki and Gorgui \cite{BG2021}, one can note that the different regimes disappear in Theorem \ref{thm:mdp-fln}. Moreover, we are able here to give the fluctuations if $2\alpha^{2} > 1$ which is not the case in \cite{BG2021}. 
\end{rem}

\section{Proof of Lemma \ref{lem:str-cv-mu}, Theorems \ref{thm:m-est-main} and \ref{thm:est-main2}  and Corollary \ref{cor:m-est-mult}}\label{sec:proof-main1}

We will  denote by  $C$ any unimportant finite  constant which may  vary from  line to line  (in particular $C$ does not  depend on $n\in \N$).

\subsection{Proof of Lemma \ref{lem:str-cv-mu}}\label{sec:Pstr-cv-mu}
We begin the proof with $\A_{n} = \TT_{n}.$ Recall the decomposition \eqref{eq:DeBiVa} with $\TT_{n}$ instead of $\A$. Using Lemma \ref{lem:bochner}, we have $\lim_{n \rightarrow +\infty} |B_{h_{n}}(x)| = 0.$ From Remark \ref{rem:cv-det-exp}, this implies that $B_{h_{n}}(x) \xRightarrow[b_{n}^{2}]{\rm superexp} 0.$ Next, we set $f_{n}(\cdot) = K_{h_{n}}(x - \cdot)$ in such a way that we have
\begin{equation*}
|\TT_{n}|^{-1} h_{n}^{d-/2} \sum_{u \in \TT_{n}} \Big(K_{h}(x - X_{u}) - K_{h}\star\mu(x)\Big) = |\TT_{n}|^{-1} h_{n}^{d-/2} \sum_{\ell = 0}^{n} M_{\GG_{\ell}}(\tilde{f}_{n}).
\end{equation*}
Following line by line the proof of \eqref{eq:I-R0n-CL} (where we take $f_{\ell,n} = f_{n}$ for all $\ell \leq n$), we get
\begin{equation*}
\PP\Big( |\TT_{n}|^{-1} h_{n}^{d-/2} \Big|\sum_{\ell = 0}^{n} M_{\GG_{\ell}}(\tilde{f}_{n})\Big| > \delta \Big) \leq 2 \exp\Big( \frac{3\delta}{c_{1} + c_{2} \delta} \Big) \exp\Big( - \frac{3 \delta^{2} |\TT_{n}| h_{n}^{d}}{c_{1} + c_{2} \delta} \Big).
\end{equation*}
Taking the $\log$, dividing by $b_{n}^{2}$ and letting $n$ goes to the infinity in the latter inequality, we get
\begin{equation*}
|\TT_{n}|^{-1} h_{n}^{d-/2} \sum_{u \in \TT_{n}} \Big(K_{h}(x - X_{u}) - K_{h}\star\mu(x)\Big) \xRightarrow[b_{n}^{2}]{\rm superexp} 0.
\end{equation*}
It then follows from the decomposition \eqref{eq:DeBiVa} that $\widehat{\mu}_{\TT_{n}}(x) \xRightarrow[b_{n}^{2}]{\rm superexp} \mu(x).$ We similarly get the result for $\A_{n} = \GG_{n}$ and this ends the proof of the lemma.

\subsection{Proof of Theorem \ref{thm:m-est-main}}\label{sec:Pm-est-main}
We begin the proof with $\A_{n} = \TT_{n}$. We have the following decomposition:
\begin{equation*}
b_{n}^{-1} \sqrt{|\TT_{n}| h_{n}^{d}}\big(\widehat{\mu}_{\TT_{n}}(x) - \mu(x)\big) = \sqrt{\frac{|\GG_{n}|}{|\TT_{n} |}} b_{n}^{-1} N_{n,\emptyset}(\bF_n) + \frac{\sqrt{|\TT_{n}| h_{n}^{d}}}{b_{n}}B_{h_n}(x),  
\end{equation*}
where $\bF_n=(f_{\ell, n}, \, n \geq \ell \geq 0)$ with the functions $f_{\ell,n}= f_{\ell,n}^{\text{id}}$ defined in \eqref{eq:def-f-kernel} for $n\geq  \ell\geq 0$ and $f_{\ell,n}=0$ otherwise; $N_{n,\emptyset}(\bF_{n})$ is defined in \eqref{eq:Nemptyf} and  the bias term $B_{h_{n}}(x)$ is defined in \eqref{eq:DeBiVa}. Thanks to Theorem \ref{thm:mdp-fln} applied to the sequence $(f_{\ell,n}^{id}, n \geq \ell \geq 0)$ and using that $\lim_{n \rightarrow +\infty} |\GG_{n}|/|\TT_{n}| = 1/2,$ we get that $\sqrt{|\GG_{n}||\TT_{n} |^{-1}} b_{n}^{-1} N_{n,\emptyset}(\bF_n)$ satisfies a moderate deviation principle in $\RR$ with speed $b_{n}^{2}$ and rate function $I$ defined by: $I(y) = y^{2}/(2 \|K\|_{2}^{2} \, \mu(x))$ for all $y \in \RR.$ To complete the proof of Theorem \ref{thm:m-est-main}, it suffices to prove that
\begin{equation}\label{eq:cv-biais-mdp}
\lim_{n \rightarrow +\infty} \frac{\sqrt{|\TT_{n}| h_{n}^{d}}}{b_{n}}B_{h_n}(x) = 0.
\end{equation}  
Next, using that 
\begin{multline*}
\mu(x - h_{n}y) - \mu(x) = \sum_{j = 1}^{d} (\mu(x_{1} - h_{n}y_{1},
  \ldots,
x_{j} - h_{n}y_{j}, x_{j+1}, \ldots, x_{d}) \\
- \mu(x_{1}-h_{n}y_{1}, \ldots, x_{j-1} - h_{n}y_{j-1}, x_{j},
    x_{j+1}, \ldots, x_{d})),  
\end{multline*}
the Taylor expansion and  Assumption \ref{hyp:estim-tcl}, 
we get that, for some finite constant $C > 0$,
\begin{align*}
|\TT_{n}|^{1/2}h_{n}^{d/2} B_{h_n}(x) & = \sqrt{|\TT_{n}|h_{n}^{d}} \,\, \Big|\int_{\RR^{d}} h_{n}^{-d} K(h_n^{-1}(x-y)) \mu(y)dy - \mu(x)\Big|\\
& = \sqrt{|\TT_{n}|h_{n}^{d}} \,\,  \Big|\int_{\RR^{d}} K(y)(\mu(x - h_{n}y) - \mu(x)) \, dy\Big| \\
&\leq C \sqrt{|\TT_{n}|h_{n}^{d}}\,\, \sum_{j = 1}^{d} \,\,\int_{\RR^{d}} K(y)\frac{(h_{n}|y_{j}|)^{s}}{\lfloor s \rfloor!} dy  
\\
&\leq C \sqrt{|\TT_{n}|h_{n}^{2s + d}}. 
\end{align*}
Now, \eqref{eq:cv-biais-mdp} follows using the latter inequality and \eqref{eq:speed-mdp}. This ends the proof of Theorem \ref{thm:m-est-main} for $\A_{n} = \TT_{n}$. The proof is similar for $\A_{n} = \GG_{n}$ using $f_{\ell,n} = f^{0}_{\ell,n}.$ 

\subsection{Proof of Theorem \ref{thm:est-main2}}\label{sec:Pest-main2}.
We begin the proof with $\A_{n} = \TT_{n}$. We have the following decomposition:
\begin{equation}\label{eq:D-muTTn}
\frac{b_{n}^{-1}  \sqrt{|\TT_{n}|h_{n}^{d}} (\widehat{\mu}_{\TT_{n}}(x) - \mu(x))}{\|K\|_{2} \sqrt{\widehat{\mu}_{\A_{n}^{*}}(x)} \vee \varpi_{n}} = T_{1}(n) + T_{2}(n)
\end{equation}
where
\begin{align*}
&T_{1}(n) = (\|K\|_{2} \sqrt{\mu(x)} b_{n})^{-1} \sqrt{|\TT_{n}| h_{n}^{d}} \Big(\widehat{\mu}_{\TT_{n}}(x) - \mu(x)\Big); \\ 
&T_{2}(n) = \Big(\frac{1}{\|K\|_{2} \sqrt{\widehat{\mu}_{\A_{n}^{*}}(x)} \vee \varpi_{n}} - \frac{1}{\|K\|_{2} \sqrt{\mu(x)}}\Big) b_{n}^{-1} \sqrt{|\TT_{n}| h_{n}^{d}} \Big(\widehat{\mu}_{\TT_{n}}(x) - \mu(x)\Big).
\end{align*}

\medskip

First, we prove that 
\begin{equation}\label{eq:cv-expT2n}
T_{2}(n) \xRightarrow[b_{n}^{2}]{\rm superexp} 0.
\end{equation}
Let $\delta > 0.$ For all $r > 0$, we have
\begin{multline*}
\PP\big(|T_{2}(n)| > \delta\big) \leq \PP\big( \big|  b_{n}^{-1} \sqrt{|\TT_{n}| h_{n}^{d}} \big(\widehat{\mu}_{\TT_{n}}(x) - \mu(x)\big) \big| > \delta/r \big) \\ 
+ \PP\big( \big| \frac{1}{\|K\|_{2} \sqrt{\widehat{\mu}_{\A_{n}^{*}}(x)} \vee \varpi_{n}} - \frac{1}{\|K\|_{2} \sqrt{\mu(x)}} \big| > r \big).
\end{multline*}
This implies that (see for e.g \cite{DZ1998}, Lemma 1.2.15)
\begin{multline}\label{eq:est-main1}
\limsup_{n \rightarrow + \infty} \frac{1}{b_{n}^{2}} \log \PP\big(|T_{2}(n)| > \delta\big) \leq \max\Big\{ \limsup_{n \rightarrow +\infty} \frac{1}{b_{n}^{2}} \log \PP\big( \big|  b_{n}^{-1} \sqrt{|\TT_{n}| h_{n}^{d}} \big(\widehat{\mu}_{\TT_{n}}(x) - \mu(x)\big) \big| > \delta/r \big); \\
\limsup_{n \rightarrow + \infty} \frac{1}{b_{n}^{2}} \log \PP\big( \big| \frac{1}{\|K\|_{2} \sqrt{\widehat{\mu}_{\A_{n}^{*}}(x)} \vee \varpi_{n}} - \frac{1}{\|K\|_{2} \sqrt{\mu(x)}} \big| > r \big)\Big\}. 
\end{multline}
Using Theorem \ref{thm:m-est-main} and the contraction principle, we have
\begin{equation}\label{eq:est-main2}
\limsup_{n \rightarrow +\infty} \frac{1}{b_{n}^{2}} \log \PP\big( \big|  b_{n}^{-1} \sqrt{|\TT_{n}| h_{n}^{d}} \big(\widehat{\mu}_{\TT_{n}}(x) - \mu(x)\big) \big| > \delta/r \big) = - \frac{\delta^{2}}{2 \|K\|_{2} \mu(x) r^{2}}.
\end{equation}
Following the step 1 of the proof of Theorem 6 in \cite{BO2018} and using Lemma \ref{lem:str-cv-mu}, we can prove that
\begin{equation*}
\|K\|_{2}^{2} \, \widehat{\mu}_{\A_{n}^{*}}(x) \vee \varpi_{n}^{2} \xRightarrow[b_{n}^{2}]{\rm superexp} \|K\|_{2}^{2} \, \mu(x).
\end{equation*}
Using Lemma B.2 in \cite{BDG14}, the latter convergence implies that
\begin{equation}\label{eq:est-main3}
\frac{1}{\|K\|_{2} \sqrt{\widehat{\mu}_{\A_{n}^{*}}(x)} \vee \varpi_{n}}  \xRightarrow[b_{n}^{2}]{\rm superexp} \frac{1}{\|K\|_{2} \sqrt{\mu(x)}}.
\end{equation}
Using \eqref{eq:est-main1}, \eqref{eq:est-main2} and \eqref{eq:est-main3}, we get
\begin{equation*}
\limsup_{n \rightarrow + \infty} \frac{1}{b_{n}^{2}} \log \PP\big(|T_{2}(n)| > \delta\big) \leq - \frac{\delta^{2}}{2 \|K\|_{2} \mu(x) r^{2}}.
\end{equation*}
Since $r$ can be taken arbitrarily close to $0$, we get \eqref{eq:cv-expT2n} and using \eqref{eq:D-muTTn}, this implies that
\begin{equation}\label{eq:est-main4}
\frac{b_{n}^{-1}  \sqrt{|\TT_{n}|h_{n}^{d}} (\widehat{\mu}_{\TT_{n}}(x) - \mu(x))}{\|K\|_{2} \sqrt{\widehat{\mu}_{\A_{n}^{*}}(x)} \vee \varpi_{n}}  \widesim[2]{{\rm superexp}}{b_{n}^{2}}  T_{1}(n).
\end{equation}
Using Theorem \ref{thm:m-est-main} and the contraction principle, we get that $T_{1}(n)$ satisfies a moderate deviation principle on $\RR$ with speed $b_{n}^{2}$ and rate function $I'$ defined by: $I'(y) = y^{2}/2$ for all $y \in \RR$. Using \eqref{eq:est-main4} and Remark \ref{rem:cv-det-exp}, we get the result of Theorem \ref{thm:est-main2}.  

\subsection{Proof of Corollary \ref{cor:m-est-mult}}\label{sec:est-mult}
Let $\boldsymbol{a} = (a_{0}, \ldots, a_{k})^{t} \in \RR^{k+1}$. Let $n > k.$ We consider the sequence $\bF_{n} = (f_{\ell,n}, n \geq \ell \geq 0)$ defined by  $f_{\ell,n} = 2^{\ell/2} \, a_{\ell} \, K_{h_{n-\ell}}(x - \cdot)$ for all $\ell \in \{0, \ldots, k\}$ and $f_{\ell,n} = 0$ otherwise. We easily check that $\bF_{n}$ satisfies Assumptions \ref{ass:fl-n}. In particular, the asymptotic variance defined in \eqref{eq:limf-ln} is given by 
$\sigma^{2} = \big(\sum_{\ell = 0}^{k} 2^{\ell} a_{\ell}^{2}\big) \|K\|_{2}^{2} \, \mu(x).$
Observe that the linear combinaison $M_{n}(\boldsymbol{a})$, with coefficients $\boldsymbol{a} = (a_{0}, \ldots, a_{k})^{t} \in \RR^{k+1},$ of the estimators $|\G_{n-\ell}|^{1/2} h_{n-\ell}^{d/2} (\widehat{\mu}_{\G_{n-\ell}}(x) - \mu(x))$, $\ell \in \{0, \ldots, k\}$ has the following decomposition:
\begin{equation}\label{eq:D-Mn-a}
M_{n}(\boldsymbol{a}) = N_{n,\emptyset}(\bF_{n}) + \sum_{\ell = 0}^{k} a_{\ell} \big( |\GG_{n-\ell}| \, h_{n-\ell}^{d} \big)^{1/2} B_{h_{n-\ell}}(x),
\end{equation}
where $N_{n,\emptyset}(\bF_{n})$ is defined in \eqref{eq:Nemptyf} and the $B_{h_{n-\ell}}(x)$, $\ell \in \{0, \ldots, k\}$, are defined in \eqref{eq:DeBiVa}. Applying Theorem \ref{thm:mdp-fln}, we get that  $b_{n}^{-1} \, N_{n,\emptyset}(\bF_{n})$ satisfies a moderate deviation principle on $\RR$ with speed $b_{n}^{2}$ and rate function $I_{x, \boldsymbol{a}}: \RR \rightarrow \RR$ defined by
\begin{equation}\label{eq:rate-Mn-a}
I_{x, \boldsymbol{a}}(y) = \frac{y^{2}}{2 \, \big(\sum_{\ell = 0}^{k} 2^{\ell} a_{\ell}^{2}\big) \|K\|_{2}^{2} \, \mu(x)}, \quad y \in \RR.
\end{equation}
Using \eqref{eq:cv-biais-mdp}, we have that
\begin{equation*}
\lim_{n \rightarrow +\infty} \frac{1}{b_{n}} \sum_{\ell = 0}^{k} a_{\ell} \big( |\GG_{n-\ell}| \, h_{n-\ell}^{d} \big)^{1/2} B_{h_{n-\ell}}(x) = 0.
\end{equation*}
Using Remark \ref{rem:cv-det-exp}, this implies that
\begin{equation}\label{eq:equi-mult}
\frac{1}{b_{n}} \sum_{\ell = 0}^{k} a_{\ell} \big( |\GG_{n-\ell}| \, h_{n-\ell}^{d} \big)^{1/2} B_{h_{n-\ell}}(x) \superexp 0.
\end{equation}
Using \eqref{eq:D-Mn-a} and \eqref{eq:equi-mult} we get that $b_{n}^{-1} M_{n}(\boldsymbol{a})$ and $b_{n}^{-1} N_{n,\emptyset}(\bF_{n})$ satisfy the same moderate deviation principle.  We then conclude that $b_{n}^{-1} M_{n}(\boldsymbol{a})$ satisfies a moderate deviation principle on $\RR$ with speed $b_{n}^{2}$ and rate function $I_{x, \boldsymbol{a}}$ defined in \eqref{eq:rate-Mn-a}. Since this is true for all vector $\boldsymbol{a} \in \RR^{k+1}$, that is for all the linear combinaisons of the estimators $|\G_{n-\ell}|^{1/2} h_{n-\ell}^{d/2} (\widehat{\mu}_{\G_{n-\ell}}(x) - \mu(x))$, $\ell \in \{0, \ldots, k\}$, we get the result of Corollary \ref{cor:m-est-mult}. 

\section{Proof of Theorem \ref{thm:mdp-fln}}\label{sec:proof-main2}

We begin with some notations. We will  denote by  $C$ any unimportant finite  constant which may  vary from  line to line  (in particular $C$ does not  depend on $n\in \N$ nor on the considered sequence of
functions $\bF_{n} = (f_{\ell,n}, n \geq \ell \geq 0)$). Let $(p_n, n\in \N)$ be a non-decreasing sequence of elements of $\N^*$ such that
\begin{equation*}
\lim_{n \rightarrow +\infty }  p_{n}^{3} \, b_{n}^{2} \, |\GG_{n-p_{n}}|^{-1}  = 0.
\end{equation*}
When there is no ambiguity, we write $p$ for $p_n$. 
\medskip

Let $i,j\in \T$. We write $i\preccurlyeq  j$ if $j\in i\T$. We denote by $i\wedge j$  the most recent  common ancestor of  $i$ and $j$,  which is defined  as   the  only   $u\in  \T$   such  that   if  $v\in   \T$  and $v\preccurlyeq i$, $v \preccurlyeq j$  then $v \preccurlyeq u$. We also define the lexicographic order $i\leq j$ if either $i \preccurlyeq j$ or $v0  \preccurlyeq i$  and $v1  \preccurlyeq j$  for $v=i\wedge  j$.  Let $X=(X_i, i\in  \T)$ be  a $BMC$  with kernel  $\cp$ and  initial measure $\nu$. For $i\in \T$, we define the $\sigma$-field:
\begin{equation*}\label{eq:field-Fi}
\cf_{i}=\{X_u; u\in \T \text{ such that  $u\leq i$}\}.
\end{equation*}
By construction,  the $\sigma$-fields $(\cf_{i}; \, i\in \T)$ are nested as $\cf_{i}\subset \cf_{j} $ for $i\leq  j$.
\medskip

We define for $n\in \N$, $i\in \G_{n-p_n}$ and $\bF_{n}$ the martingale increments:
\begin{equation}\label{eq:def-DiF}
\Delta_{n,i}(\bF_{n})= N_{n,i}(\bF_{n}) - \E\left[N_{n,i}(\bF_{n}) |\, \cf_i\right] \quad \text{and} \quad \Delta_n(\bF_{n}) = \sum_{i\in \G_{n-p_n}} \Delta_{n,i}(\bF_{n}),
\end{equation}
where
\begin{equation}\label{eq:def-NiFn}
N_{n,i}(\bF_{n}) = |\GG_{n}|^{-1/2} \sum_{\ell = 0}^{p} M_{i\GG_{p-\ell}}(\tilde{f}_{\ell,n}) \quad \text{and} \quad i\GG_{p-\ell} = \{ij, j \in \GG_{p-\ell}\}.
\end{equation}
We have:
\begin{equation*}
\sum_{i\in \G_{n-p_n}} N_{n, i}(\bF_{n}) = |\G_n|^{-1/2} \sum_{\ell=0}^{p_n}  M_{\G_{n-\ell}} (\tilde f_{\ell,n}) = |\G_n|^{-1/2} \sum_{k=n-p_n}^{n}  M_{\G_{k}} (\tilde f_{n-k,n}).
\end{equation*}
Using the branching Markov property, we get for $i\in \G_{n-p_n}$:
\begin{equation}\label{eq:EN-ni-fn}
\E\left[N_{n,i}(\bF_{n}) |\, \cf_i\right] =\E\left[N_{n,i}(\bF_{n}) |\, X_i\right] = |\G_n|^{-1/2} \sum_{\ell=0}^{p_n} \E_{X_i}\left[M_{\G_{p_n-\ell}}(\tilde f_{\ell,n})\right].
\end{equation}
We have the following decomposition:
\begin{equation}\label{eq:D-N0fn}
N_{n, \emptyset}(\bF_{n}) = \Delta_n(\bF_{n}) + R_{0} (n)+R_{1}(n),
\end{equation}
where $\Delta_n(\bF)$ is defined in \reff{eq:def-DiF} and:
\begin{equation*}
R_{0} (n) = |\G_n|^{-1/2}\, \sum_{k = 0}^{n-p_n-1} M_{\G_k}(\tilde {f}_{n-k,n}) \quad \text{and} \quad R_{1}(n) = \sum_{i\in \G_{n-p_n}}\E\left[N_{n,i}(\bF_{n}) |\, \cf_i\right].
\end{equation*}
From \eqref{eq:D-N0fn}, our goals will be achieved if we prove the following:
\begin{align}
&\label{eq:cvsupexpR0}
b_{n}^{-1} R_{0}(n) \superexp 0;\\
&\label{eq:cvsupexpR1}
b_{n}^{-1} R_{1}(n) \superexp 0;\\
&\label{eq:mdp-Dnf-SCr}
b_{n}^{-1}\Delta_{n}(\bF) \quad \text{satisfies a MDP on $S$ with speed $b_n^2$ and rate function $I$.}
\end{align}
Note that \eqref{eq:cvsupexpR0} and \eqref{eq:cvsupexpR1} mean that $R_{0}(n)$ and $R_{1}(n)$ are negligible in the sense of moderate deviations in such a way that using \eqref{eq:D-N0fn} and Remark \ref{rem:cv-det-exp}, $N_{n,\emptyset}(\bF)$ and $\Delta_{n}(\bF)$ satisfy the same moderate deviation principle. To prove \eqref{eq:mdp-Dnf-SCr}, the main method we will use is the moderate deviations for martingale (see \cite{djellout2002moderate} for more details).

\medskip

In the sequel, the sequence $(2^{-\gamma n}, n \in \NN)$ which appears in Assumption \ref{ass:fl-n} will be denoted $(h_{n}, n \in \NN)$ in such a way that we have $2^{-d \gamma n/2} = h_{n}^{d/2}.$ We have the following result.
\begin{lem}\label{lem:cvexp-R0n}
Under the assumptions of Theorem \ref{thm:mdp-fln}, we have $b_{n}^{-1} R_{0}(n) \superexp 0.$
\end{lem}
\begin{proof}
Let $\delta > 0$. Using the Chernoff bound, we have, for all $\lambda > 0$,
\begin{equation}\label{eq:BC-R0n}
\PP\Big(b_{n}^{-1} R_{0}(n) > \delta\Big) \leq \exp\Big(- \lambda b_{n} |\GG_{n}|^{1/2} \delta\Big) \EE\Big[\exp\Big(\lambda \sum_{\ell = 0}^{n-p-1} M_{\GG_{\ell}}(\tilde{f}_{n-\ell,n})\Big)\Big].
\end{equation}
For all $k \in \{1,\ldots,n-p\}$ and for $u \in \TT$, we set
\[ 
g_{p,k} = \sum_{r = 0}^{k-1} 2^{r} \Qq^{r} \tilde{f}_{p+k-r,n} \quad \text{and} \quad Z_{p,k}(u) = g_{p,k}(X_{u0}) + g_{p,k}(X_{u1}) - 2 \Qq g_{p,k}(X_{u}).
\]
Then, using recursively the fact that 
\begin{equation*}
\sum_{u \in \GG_{\ell}} f(X_{u}) =  \sum_{u \in \GG_{\ell -1}} (f(X_{u0}) + f(X_{u1}) - 2\cq f(X_{u})) + \sum_{u \in \GG_{\ell - 1}} 2\cq f(X_{u}),
\end{equation*}
for all $\ell \geq 1$ and for some function $f$, we get
\begin{equation*}
\EE\Big[\exp\Big(\lambda \sum_{\ell = 0}^{n-p-1} M_{\GG_{\ell}}(\tilde{f}_{n-\ell,n})\Big)\Big] = \EE\Big[\exp\Big(\lambda g_{p,n-p}(X_{\emptyset})\Big) \, \prod_{k=1}^{n-p-1} \exp\Big(\lambda \sum_{u \in \GG_{n-p-k-1}} Z_{p,k}(u)\Big)\Big].
\end{equation*}
For all $m \in \{1, \ldots, n-p-1\}$,  we set
\begin{equation*}
\II_{m} = \EE\Big[\exp(\lambda g_{p,n-p}(X_{\emptyset}))  \, \prod_{k = m}^{n-p-1} \exp(\lambda \sum_{u \in \GG_{n-p-k-1}} Z_{p,k}(u))\Big].
\end{equation*}
Using the branching Markov property, we get the following decomposition:
\begin{equation*}
\II_{m} = \EE\Big[\exp(\lambda g_{p,n-p}(X_{\emptyset})) \, \JJ_{m} \, \prod_{k = m + 1}^{n-p-1} \exp(\lambda \sum_{u \in \GG_{n-p-k-1}} Z_{p,k}(u))\Big],
\end{equation*}
with
\begin{equation*}
\JJ_{m} = \prod_{u \in \GG_{n-p-m-1}}\EE_{X_{u}}\Big[\exp(\lambda Z_{p,m}(u))\Big].
\end{equation*}
For all $u \in \GG_{n-p-m-1}$, we will upper bound the quantity $\EE_{X_{u}}[\exp(\lambda Z_{p,m}(u))]$ and then $\JJ_{m}$. We claim that:
\begin{equation}\label{eq:claim1-R0n}
|Z_{p,m}(u)| \, \leq \, M \, = \, C \, h_{n}^{-d/2};
\end{equation}
\begin{equation}\label{eq:claim2-R0n}
\EE_{X_{u}}[Z_{p,m}(u)^{2}] \, \leq \, \sigma_{m}^{2} \, = \,C +  \, C \, h_{n}^{d} \, \Big(\sum_{r = 0}^{m-1} (2\alpha)^{r-1}\Big)^{2} \ind_{\{m>1\}}.
\end{equation}

For that purpose, we plan to use the bound
\begin{equation} \label{eq:Bennett}
\EE\big[\exp(\lambda Z) \big] \leq \exp\Big(\frac{\lambda^2 \sigma^2}{2(1 - \lambda M / 3)}\Big)
\end{equation}
valid for any $\lambda \in(0,3/M)$, any  random variable $Z$ such that $|Z| \leq M$, $\EE[Z] = 0$ and $\EE[Z^2] \leq \sigma^2$. For all $u \in \GG_{n-p-m-1}$ and for all $\lambda \in (0,C h_{n}^{-d/2}/3)$ we get, using \eqref{eq:claim1-R0n}-\eqref{eq:Bennett},
\begin{equation*}
\EE_{X_{u}}\Big[\exp(\lambda Z_{p,m}(u))\Big] \leq \exp\Big(\frac{\lambda^{2} \sigma_{m}^{2}}{2(1 -  \lambda M/3)}\Big).
\end{equation*}
For all $m \in \{1, \ldots, n-p-1\}$, the latter inequality implies that
\begin{equation*}
\JJ_{m} \leq \exp\Big(\frac{\lambda^{2} \sigma_{m}^{2} |\GG_{n-p-m-1}|}{2(1 - \lambda M/3)}\Big) \quad \text{and} \quad \II_{m} \leq \exp\Big(\frac{\lambda^{2} \sigma_{m}^{2} |\GG_{n-p-m-1}|}{2(1 - \lambda M/3)}\Big) \II_{m+1}.
\end{equation*}
Recall that $\II_{1} = \EE\big[\exp\big(\lambda \sum_{\ell = 0}^{n-p-1} M_{\GG_{\ell}}(\tilde{f}_{n-\ell,n})\big)\big].$  By recurrence, we get
\begin{equation*}
\EE\Big[\exp\Big(\lambda \sum_{\ell = 0}^{n-p-1} M_{\GG_{\ell}}(\tilde{f}_{n-\ell,n})\Big)\Big] = \II_{1} \leq \exp\Big(\frac{\lambda^{2} \sum_{m = 1}^{n - p - 1} \sigma_{m}^{2} |\GG_{n-p-m-1}|}{2(1 - \lambda M/3)}\Big) \EE\Big[\exp\Big(\lambda g_{p,n-p}(X_{\emptyset})\Big)\Big].
\end{equation*}
Using $(i)$ and $(ii)$ of Assumption \ref{ass:fl-n} and \eqref{eq:geom-erg}, we have
\begin{align*}
|g_{p,n-p}| \leq |\tilde{f}_{n,n}| + \sum_{r = 1}^{n-p-1} 2^{r} |\Qq^{r-1} (\Qq\tilde{f}_{n-r,n})| \leq C h_{n}^{-d/2} + C h_{n}^{d/2} \sum_{r = 1}^{n-p-1} (2\alpha)^{r-1}.
\end{align*}
This implies that
\begin{multline*}
\EE\Big[\exp\Big(\lambda \sum_{\ell = 0}^{n-p-1} M_{\GG_{\ell}}(\tilde{f}_{n-\ell,n})\Big)\Big] \leq \exp\Big(\frac{\lambda^{2} \sum_{m = 1}^{n - p - 1} \sigma_{m}^{2} |\GG_{n-p-m-1}|}{2(1 - \lambda M/3)}\Big) \\ \times \, \exp\Big(\lambda \, C \, h_{n}^{-d/2} \,  + \, \lambda \, C \,  h_{n}^{d/2} \, \sum_{r = 0}^{n-p-2} (2\alpha)^{r}\Big).
\end{multline*}
Distinguishing the cases $2\alpha \leq 1$, $1/2 < 2\alpha \leq \sqrt{2}$ and $2\alpha > \sqrt{2}$ and using \eqref{eq:fln-aS2} for $2\alpha > 1$, we get
\begin{equation*}
\EE\Big[\exp\Big(\lambda \sum_{\ell = 0}^{n-p-1} M_{\GG_{\ell}}(\tilde{f}_{n-\ell,n})\Big)\Big] \leq \exp\Big( \frac{c_{1} \lambda^{2} |\GG_{n-p}|}{2(1 - c_{2} \lambda h_{n}^{-d/2}/3)} \Big) \exp\Big( c_{3} \lambda h_{n}^{-d/2} \Big),
\end{equation*}
where $c_{1}$, $c_{2}$ and $c_{3}$ are some positive constants. The latter inequality and \eqref{eq:BC-R0n} imply that
\begin{equation*}
\PP\Big(b_{n}^{-1} R_{0}(n) > \delta\Big) \leq \exp\Big(- \lambda b_{n} |\GG_{n}|^{1/2} \delta + \frac{c_{1} \lambda^{2} |\GG_{n-p}|}{2(1 - c_{2} \lambda h_{n}^{-d/2}/3)} \Big) \, \exp\Big( c_{3} \lambda h_{n}^{-d/2} \Big).
\end{equation*}
Taking\footnote{In fact, we use the following. For $\alpha, \beta, \gamma > 0$ and $h(x)=- \alpha x + \frac{\beta x^2}{2 (1 - \gamma x)}$ we have $h(x^*) = \frac{-\alpha^2}{2(\beta + \alpha \gamma)}$ for the choice $x^* = \frac{\alpha}{2 \alpha \gamma + \beta} \in (0,1/\gamma).$}
\begin{equation*}
\lambda = \frac{3 \, b_{n} \, |\GG_{n}|^{1/2} \, \delta}{2 \, c_{2} \, b_{n} \, |\GG_{n}|^{1/2} \, h_{n}^{-d/2} \,\delta + 3 \, c_{1} \, |\GG_{n-p}|}, 
\end{equation*}
we are led to
\begin{equation*}
\PP\Big(b_{n}^{-1} R_{0}(n) > \delta\Big) \leq C \, \exp\Big( - \,\frac{3 \, \delta^{2} \, b_{n}^{2} \, |\GG_{n}|}{2( c_{2} \, \delta \, b_{n} \, |\GG_{n}|^{1/2} \, h_{n}^{-d/2} + 3 \, c_{1} \, |\GG_{n-p}|)} \Big).
\end{equation*}
Since we can do the same thing for $-\bF_{n}$ instead of $\bF_{n}$, we get that
\begin{equation}\label{eq:I-R0n-CL}
\PP\Big(b_{n}^{-1} |R_{0}(n)| > \delta\Big) \leq 2 \, C \, \exp\Big( - \,\frac{3 \, \delta^{2} \, b_{n}^{2} \, |\GG_{n}|}{2( c_{2} \, \delta \, b_{n} \, |\GG_{n}|^{1/2} \, h_{n}^{-d/2} + 3 \, c_{1} \, |\GG_{n-p}|)} \Big).
\end{equation}
Finally, in the latter inequality, taking the $\log$, dividing by $b_{n}^{2}$ and letting $n$ goes to infinity, we get the result of Lemma \ref{lem:cvexp-R0n}. Now, to end the proof, we will prove \eqref{eq:claim1-R0n} and \eqref{eq:claim2-R0n}.

\subsubsection*{Proof of \eqref{eq:claim1-R0n}} Using Assumption \ref{hyp:erg-unif}, $(i)$ and $(ii)$ of Assumption \ref{ass:fl-n} and Assumption \ref{ass:fln-aS2}, we get
\begin{align*}
|Z_{p,m}(u)| &\leq C \, \|\tilde{f}_{p+1,n}\|_{\infty}  +  C(1+2\alpha)( \sum_{r = 1}^{m-1} (2\alpha)^{r-1} \|\Qq f_{p+m-r,n}\|_{\infty}) \ind_{\{m > 1\}} \\  &\leq C \, h_{n}^{-d/2} \, + \, C \, h_{n}^{d/2} \, \sum_{r = 0}^{m-1} (2\alpha)^{r} \, \leq \, C \, h_{n}^{-d/2} .
\end{align*}
\subsubsection*{Proof of \eqref{eq:claim2-R0n}} Using the branching Markov property for the second inequality, Assumption \ref{hyp:erg-unif} for the fourth inequality and $(i)$ and $(ii)$ of Assumption \ref{ass:fl-n} for the last inequality, we get
\begin{align*}
\EE_{X_{u}}[Z_{p,m}(u)^{2}] &\leq \EE_{X_{u}}[(g_{p,m}(X_{u0}) + g_{p,m}(X_{u1}))^{2}] \leq  C \Qq(g_{p,m}^{2})(X_{u}) \\ & \leq C \Qq(\tilde{f}_{p+1,n}^{2})(X_{u}) + C \Qq\Big(\Big(\sum_{r = 1}^{m-1} 2^{r} \Qq^{r-1} (\Qq \tilde{f}_{p+m-r,n})\Big)^{2}\Big)(X_{u}) \, \ind_{\{m > 1\}} \\ & \leq C \|\Qq \tilde{f}_{p+1,n}^{2}\|_{\infty} + \Big( \sum_{r=1}^{m-1} (2\alpha)^{r-1} \|\Qq f_{p+m-r,n}\|_{\infty} \Big)^{2} \ind_{\{m>1\}} \\ & \leq \, C \, + \, C \, h_{n}^{d} \, \Big(\sum_{r=0}^{m-1} (2\alpha)^{r}\Big)^{2} \ind_{\{m>1\}}.
\end{align*}
\end{proof}

Next, we have the following result.
\begin{lem}\label{lem:cvexp-R1n}
Under the assumptions of Theorem \ref{thm:mdp-fln}, we have $b_{n}^{-1} R_{1}(n) \superexp 0.$
\end{lem}
\begin{proof}
We have, using \eqref{eq:EN-ni-fn} and \eqref{eq:Q1},
\begin{equation}\label{eq:R1n-gpn}
R_{1}(n) = |\GG_{n}|^{-1/2} M_{\GG_{n-p}}(g_{p,n}) \quad \text{where} \quad g_{p,n} = \sum_{\ell = 0}^{p} 2^{p-\ell} \Qq^{p-\ell} \tilde{f}_{\ell,n}.
\end{equation}
We follow the same arguments that in the proof of Lemma \ref{lem:cvexp-R0n}. For all $m \in \{1 ,\ldots, n-p\}$ and for all $u \in \TT$, we set
\begin{equation*}
Z_{p,m}(u) = 2^{m-1} \Qq^{m-1} g_{p}(X_{u0}) + 2^{m-1} \Qq^{m-1} g_{p}(X_{u1}) - 2^{m} \Qq^{m} g_{p}(X_{u}).
\end{equation*}
We also consider the following quantities for $m \in \{1, \ldots, n-p\}$ and $\lambda > 0$:
\begin{align*}
&\II_{m} = \EE\Big[\exp\Big(\lambda 2^{n-p} \Qq^{n-p} g_{p,n}(X_{\emptyset})\Big)\prod_{k = m}^{n-p} \exp\Big(\lambda \sum_{u \in \GG_{n-p-k}} Z_{p,k}(u)\Big)\Big] \quad \text{and} \\ &\JJ_{m} = \prod_{u \in \GG_{n-p-m}} \EE_{X_{u}}\Big[\exp\Big(\lambda Z_{p,m}(u)\Big)\Big]. 
\end{align*}
Note that using the branching Markov property, we have
\begin{equation}\label{eq:IIm-JJm-R1n}
\II_{m} = \EE\Big[\exp\Big(\lambda 2^{n-p} \Qq^{n-p} g_{p,n}(X_{\emptyset})\Big)\prod_{k = m + 1}^{n-p} \exp\Big(\lambda \sum_{u \in \GG_{n-p-k}} Z_{p,k}(u)\Big) \JJ_{m}\Big].
\end{equation}
As for \eqref{eq:claim1-R0n}-\eqref{eq:claim2-R0n}, for all $m \in \{1, \ldots, n-p\}$ and $u \in \GG_{n-p-m}$, one can prove that
\begin{equation}\label{eq:claim1-R1n}
|Z_{p,m}(u)| \leq M = C h^{-d/2} \quad \text{and} \quad \EE_{X_{u}}\Big[Z_{p,m}(u)^{2}\Big] \leq \sigma_{m}^{2} = C \ind_{\{m=1\}} + C h_{n}^{d} \Big(\sum_{\ell = 0}^{p} (2\alpha)^{p+m-\ell-2}\Big)^{2}.
\end{equation}
Using \eqref{eq:Bennett} and \eqref{eq:claim1-R1n}, we have, for all $u \in \GG_{n-p-m}$ and for all $\lambda \in (0, C h^{-d/2}/3)$,
\begin{equation*}
\EE_{X_{u}}\Big[\exp(\lambda Z_{p,m}(u))\Big] \leq \exp\Big(\frac{\lambda^{2} \sigma_{m}^{2}}{2(1 -  \lambda M/3)}\Big).
\end{equation*}
The latter inequality and \eqref{eq:IIm-JJm-R1n} imply that
\begin{equation*}
\II_{m} \leq \exp\Big(\frac{\lambda^{2} \sigma_{m}^{2} |\GG_{n-p-m}|}{2(1 - \lambda M/3)}\Big) \II_{m+1}.
\end{equation*}
By recurrence, this implies that
\begin{equation}\label{eq:II-1}
\II_{1} \leq \exp\Big(\frac{\lambda^{2} \sum_{m=1}^{n-p} \sigma_{m}^{2} |\GG_{n-p-m}|}{2(1 - \lambda M/3)}\Big) \, \EE\Big[\exp\Big(\lambda 2^{n-p} \Qq^{n-p} g_{p,n}(X_{\emptyset})\Big)\Big].
\end{equation}
Using $(i)$ and $(ii)$ of Assumption \ref{ass:fl-n} and Assumption \ref{hyp:erg-unif}, we get
\begin{equation}\label{eq:g-p-n}
|g_{p,n}| \leq C h_{n}^{d/2} \sum_{\ell = 0}^{p} (2\alpha)^{n-\ell}.
\end{equation}
From \eqref{eq:II-1}, \eqref{eq:g-p-n} and according to the value of $\alpha$, we have, for some positive constants $c_{1}$, $c_{2}$ and $c_{3}$ (recall the definition of $M$ and $\sigma_{m}^{2}$ given in \eqref{eq:claim1-R1n}):
\begin{align*}
&\II_{1} \leq C \exp\Big( \frac{\lambda^{2} c_{1} |\GG_{n-p}|}{2(1 - \lambda c_{2} h_{n}^{-d/2}/3)} \Big) \quad \hspace{6.2cm} \text{if $2\alpha \leq 1;$}\\
&\II_{1} \leq \exp\Big(\lambda c_{3} (2\alpha)^{n}h_{n}^{d/2}\Big) \, \exp\Big( \frac{\lambda^{2} c_{1} |\GG_{n-p}|(1 + (2\alpha)^{2p}h_{n}^{d})}{2(1 - \lambda c_{2} h_{n}^{-d/2}/3)} \Big) \quad \hspace{2.2cm} \text{if $1 < 2\alpha \leq \sqrt{2};$}\\
&\II_{1} \leq \exp\Big(\lambda c_{3} (2\alpha)^{n}h_{n}^{d/2}\Big) \, \exp\Big( \frac{\lambda^{2} c_{1} |\GG_{n-p}|(1 + (2\alpha)^{2p}h_{n}^{d} + 2^{p}(2\alpha^{2})^{n}h_{n}^{d})}{2(1 - \lambda c_{2} h_{n}^{-d/2}/3)} \Big) \quad \text{if $2\alpha > \sqrt{2}.$}
\end{align*}
 
Recall that $\II_{1} = \EE[\exp(\lambda M_{\GG_{n-p}}(g_{p,n}))]$. Using the Chernoff bound and \eqref{eq:R1n-gpn}, we have for all $\lambda \in (0, C h_{n}^{-d/2}/3)$ and for all $\delta > 0$,
\begin{equation*}
\PP\Big(b_{n}^{-1} R_{1}(n) > \delta\Big) \leq \exp\Big(- \lambda b_{n} |\GG_{n}|^{1/2} \delta\Big) \,   \II_{1}.
\end{equation*}

Taking
\begin{equation*}
\lambda = \begin{cases} \frac{3  b_{n}  |\GG_{n}|^{1/2}  \delta}{2  c_{2}  b_{n}  |\GG_{n}|^{1/2}  h_{n}^{-d/2} \delta \, + \, 3  c_{1}  |\GG_{n-p}|} & \text{if $2\alpha \leq 1$} \\ \frac{3  b_{n}  |\GG_{n}|^{1/2}  \delta}{2  c_{2}  b_{n}  |\GG_{n}|^{1/2}  h_{n}^{-d/2} \delta \, + \, 3  c_{1}  |\GG_{n-p}|(1 + (2\alpha)^{2p}h_{n}^{d})} & \text{if $1 < 2\alpha \leq \sqrt{2}$} \\ \frac{3  b_{n}  |\GG_{n}|^{1/2}  \delta}{2  c_{2}  b_{n}  |\GG_{n}|^{1/2}  h_{n}^{-d/2} \delta \, + \, 3  c_{1}  |\GG_{n-p}|(1 + (2\alpha)^{2p}h_{n}^{d} + 2^{p}(2\alpha^{2})^{n}h_{n}^{d})} & \text{if $1 < 2\alpha \leq \sqrt{2},$}   \end{cases}
\end{equation*}
and since we can do the same things for $-\bF_{n}$ instead of $\bF_{n}$, we get, 

\noindent if $2\alpha \leq 1:$
\begin{equation*}
\PP\Big(b_{n}^{-1} |R_{1}(n)| > \delta\Big) \leq C \exp\Big(- \frac{3  b_{n}^{2}  |\GG_{n}|  \delta^{2}}{2  (c_{2}  b_{n}  |\GG_{n}|^{1/2}  h_{n}^{-d/2} \delta \, + \, 3  c_{1}  |\GG_{n-p}|)} \Big);
\end{equation*}
if $1 < 2\alpha\leq\sqrt{2}:$
\begin{align*}
\PP\Big(b_{n}^{-1} |R_{1}(n)| > \delta\Big) &\leq 2 \exp\Big( \frac{c_{3} (2\alpha)^{n} h_{n}^{d/2} b_{n} |\GG_{n}|^{1/2}}{2 c_{2} b_{n} |\GG_{n}|^{1/2} h_{n}^{-d/2} \delta \, +  \, 3 c_{1} |\GG_{n-p}|(1 + (2\alpha)^{2p} h_{n}^{d})} \Big) \\& \hspace{1cm} \times \exp\Big(- \frac{3  b_{n}^{2}  |\GG_{n}|  \delta^{2}}{2  (c_{2}  b_{n}  |\GG_{n}|^{1/2}  h_{n}^{-d/2} \delta \, + \, 3  c_{1}  |\GG_{n-p}|(1 + (2\alpha)^{2p}h_{n}^{d}))}\Big);
\end{align*}
if $2\alpha > \sqrt{2}:$
\begin{align*}
\PP\Big(b_{n}^{-1} |R_{1}(n)| > \delta\Big) &\leq 2 \exp\Big( \frac{c_{3} (2\alpha)^{n} h_{n}^{d/2} b_{n} |\GG_{n}|^{1/2}}{2 c_{2} b_{n} |\GG_{n}|^{1/2} h_{n}^{-d/2} \delta \, +  \, 3 c_{1} |\GG_{n-p}|(1 + (2\alpha)^{2p} h_{n}^{d} + 2^{p} (2\alpha^{2})^{n}h_{n}^{d})} \Big) \\& \hspace{0.25cm} \times \exp\Big(- \frac{3  b_{n}^{2}  |\GG_{n}|  \delta^{2}}{2  (c_{2}  b_{n}  |\GG_{n}|^{1/2}  h_{n}^{-d/2} \delta \, + \, 3  c_{1}  |\GG_{n-p}|(1 + (2\alpha)^{2p}h_{n}^{d} + 2^{p} (2\alpha^{2})^{n}h_{n}^{d}))}\Big).
\end{align*}
Finally, applying the $\log$ to each of these last three inequalities, dividing by $b_{n}^{2}$, letting $n$ goes to infinity and using \eqref{eq:speed-mdp} and Assumption \ref{ass:fln-aS2}, we get the result of Lemma \ref{lem:cvexp-R1n}.
\end{proof}

From \eqref{eq:D-N0fn}, Lemmas \ref{lem:cvexp-R0n} and \ref{lem:cvexp-R1n}, we have
\begin{equation}\label{eq:equiexp-N-D} 
b_{n}^{-1} N_{n, \emptyset}(\bF_{n}) \equiexp b_{n}^{-1} \Delta_{n}(\bF_{n}). 
\end{equation}
As a consequence, using Remark \ref{rem:cv-det-exp}, $b_{n}^{-1} N_{n, \emptyset}(\bF_{n})$ and $b_{n}^{-1} \Delta_{n}(\bF_{n})$ satisfy the same moderate deviation principle.

\medskip

We now study the martingale part $\Delta_{n}(\bF_{n})$ of the decomposition \eqref{eq:D-N0fn}. The bracket $V(n)$ of $\Delta_n(\bF_{n})$ is defined by:
\begin{equation*}\label{eq:def-Vn-subc}
V(n)= \sum_{i\in \G_{n-p_n}} \E\left[ \Delta_{n, i}(\bF_{n})^2|\cf_i\right]. 
\end{equation*}
Using \reff{eq:def-NiFn} and \reff{eq:def-DiF}, we write:
\begin{equation}\label{eq:def-V}
V(n) = |\G_n|^{-1} \sum_{i\in \G_{n-p_n}} \E_{X_i}\left[\left(\sum_{\ell=0}^{p_n} M_{\G_{p_n-\ell}}(\tilde f_{\ell,n}) \right)^2 \right]-R_2(n)=V_1(n) +2V_2(n) - R_2( n),
\end{equation}
with:
\begin{align*}
V_1(n) & =   |\G_n|^{-1} \sum_{i\in \G_{n-p_n}} \sum_{\ell=0}^{p_n} \E_{X_i}\left[M_{\G_{p_n-\ell}}(\tilde f_{\ell,n}) ^2  \right] ,\\
V_2(n) & =  |\G_n|^{-1} \sum_{i\in \G_{n-p_n}} \sum_{0\leq \ell<k\leq p_n} \E_{X_i}\left[M_{\G_{p_n-\ell}}(\tilde f_{\ell,n})  M_{\G_{p_n-k}}(\tilde f_{k,n}) \right], \\
R_2( n) &=\sum_{i\in \G_{n-p_n}} \E\left[ N_{n,i} (\bF_{n}) |X_i \right] ^2.
\end{align*}
We have the following result.
\begin{lem}\label{lem:cvR2nfln}
Under the Assumptions of Theorem \ref{thm:mdp-fln}, we have $R_{2}(n) \superexp 0.$
\end{lem}
\begin{proof}
Using the branching Markov property, we have
\begin{equation*}
R_{2}(n) = |\GG_{n}|^{-1} M_{\GG_{n-p}}(g_{p}) \quad \text{with} \quad g_{p} = \Big(\sum_{\ell = 0}^{p} 2^{p-\ell} \Qq^{p-\ell} \tilde{f}_{\ell,n}\Big)^{2}.
\end{equation*}
Using Assumption \ref{hyp:erg-unif} and $(i)$ and $(ii)$ of Assumption \ref{ass:fl-n}, we get
\begin{align*}
\|g_{p}\|_{\infty} &\leq C \|\tilde{f}_{p,n}\|^{2}_{\infty} + C \|(\sum_{\ell = 0}^{p-1} 2^{p-\ell} \Qq^{p-\ell} \tilde{f}_{\ell,n})^{2}\|_{\infty} \\
&\leq C h_{n}^{-d} + C \Big( \sum_{\ell = 0}^{p-1} (2\alpha)^{p-\ell}h_{n}^{d/2} \Big)^{2} \\
&\leq C h_{n}^{-d} \ind_{\{2\alpha \leq 1\}} + C(h_{n}^{-d} + h_{n}^{d} (2\alpha)^{2p}) \ind_{\{2\alpha > 1\}}.
\end{align*}
This implies that
\begin{equation}\label{eq:RA1-2}
R_{2}(n) \leq C |\GG_{n}|^{-1}h_{n}^{-d} \, \ind_{\{2\alpha \leq 1\}} +  C (|\GG_{n}|^{-1}h_{n}^{-d} + (2\alpha^{2})^{p} h_{n}^{d} |\GG_{n-p}|^{-1}) \, \ind_{\{2\alpha > 1\}}.
\end{equation}
Recall that $h_{n} = 2^{-n\gamma}$ with $\gamma \in (0,1/d)$. Using Assumption \ref{ass:fln-aS2}, we conclude from \eqref{eq:RA1-2} that $R_{2}(n)$ is bounded by a deterministic sequence which converge to 0. As a consequence, using Remark \ref{rem3}, we get the result of Lemma \ref{lem:cvR2nfln}.
\end{proof}
Recall $\sigma^{2}$ given in \eqref{eq:limf-ln}. We have the following result.
\begin{lem}\label{lem:cvV1nfln}
Under the Assumptions of Theorem \ref{thm:mdp-fln}, we have $V_{1}(n) \superexp \sigma^{2}.$
\end{lem}
\begin{proof}
We have the following decomposition which is a consequence of \reff{eq:Q2}:
\begin{equation*}
V_1(n)= V_3(n)+ V_4(n),
\end{equation*}
with
\begin{align*}
V_3(n) &=  |\G_n|^{-1} \sum_{i\in \G_{n-p}} \sum_{\ell=0}^p 2^{p-\ell}\,\cq^{p-\ell} (\tilde f_{\ell,n}^2)(X_i),\\ 
V_4(n) &=     |\G_n|^{-1} \sum_{i\in \G_{n-p}} \sum_{\ell=0}^{p-1}\,\sum_{k=0}^{p-\ell -1} 2^{p-\ell+k} \,\cq^{p-1-(\ell+k)}\left(\cp\left(\cq^k \tilde f_{\ell,n} \otimes^2\right)\right)(X_i).  
\end{align*}
Now, the result of Lemma \ref{lem:cvV1nfln} is a direct consequence of the following:
\begin{align}\label{eq:cvV3n-fln}
&V_{3}(n) \superexp \sigma^{2} 
;\\
&\label{eq:cvV4n-fln} V_{4}(n) \superexp 0
.
\end{align}
To end the proof, we will now prove \eqref{eq:cvV3n-fln} and \eqref{eq:cvV4n-fln}.
\subsubsection*{Proof of \eqref{eq:cvV3n-fln}}

Set
\[ 
g_{p,n} = \sum_{\ell = 0}^{p} 2^{-\ell} \Qq^{p-\ell}(\tilde{f}^{2}_{\ell,n} - \langle \mu, \tilde{f}^{2}_{\ell,n}\rangle) \quad \text{and} \quad H_{3}^{[n]}(\bF_{n}) = \sum_{\ell = 0}^{p} 2^{-\ell} \langle \mu,\tilde{f}_{\ell,n}^{2} \rangle.
\]
Following the same arguments that in the proof of Lemmas \ref{lem:cvexp-R0n} and \ref{lem:cvexp-R1n}, we get after studious calculations:

\noindent if $2\alpha \leq 1,$
\begin{align*}
&\PP\Big(|V_{3}(n) - H_{3}^{[n]}| > \delta\Big) = \PP\Big(|\GG_{n-p}|^{-1} |M_{\GG_{n-p}}(g_{p,n})| > \delta\Big) \\
& \leq C \exp\Big(\frac{C  p \, \delta}{C \delta h_{n}^{-d} + 3(p^{2} 2^{-p} + 2^{-p} h_{n}^{-d})}\Big) \exp\Big(- \frac{3\delta^{2} |\GG_{n}|}{2(C \delta h_{n}^{-d} + 3(p^{2}2^{-p} + 2^{-p}h_{n}^{-d}))}\Big);
\end{align*}
if $1 < 2\alpha \leq \sqrt{2},$
\begin{align*}
&\PP\Big(|V_{3}(n) - H_{3}^{[n]}| > \delta\Big) = \PP\Big(|\GG_{n-p}|^{-1} |M_{\GG_{n-p}}(g_{p,n})| > \delta\Big) \\
& \leq \exp\Big(\frac{c_{3} (2\alpha)^{n}h_{n}^{d} \delta}{c_{2} \delta  + 3 c_{1}((2\alpha^{2})^{p}h_{n}^{d} + 2^{-p})}\Big) \exp\Big(-  \frac{3 \delta^{2} |\GG_{n}| h_{n}^{d}}{2(c_{2} \delta + 3 c_{1}((2\alpha^{2})^{p} h_{n}^{d} + 2^{-p}))}\Big);
\end{align*}
if $2\alpha > \sqrt{2},$
\begin{align*}
&\PP\Big(|V_{3}(n) - H_{3}^{[n]}| > \delta\Big) = \PP\Big(|\GG_{n-p}|^{-1} |M_{\GG_{n-p}}(g_{p,n})| > \delta\Big) \\
& \leq \exp\Big(\frac{c_{3} (2\alpha)^{n}h_{n}^{d} \delta}{c_{2} \delta  + 3 c_{1}((2\alpha^{2})^{n}h_{n}^{d} + 2^{-p})}\Big) \exp\Big(-  \frac{3 \delta^{2} |\GG_{n}| h_{n}^{d}}{2(c_{2} \delta + 3 c_{1}((2\alpha^{2})^{n} h_{n}^{d} + 2^{-p}))}\Big);
\end{align*}
Taking the $\log$, dividing by $b_{n}^{2}$, letting $n$ goes to the infinity and using \eqref{eq:speed-mdp} and Assumption \ref{ass:fln-aS2}, we get
\begin{equation*}
\limsup_{n \rightarrow +\infty} \frac{1}{b_{n}^{2}} \log \PP\Big(|V_{3}(n) - H_{3}^{[n]}| > \delta\Big) = -\infty.
\end{equation*}
Next, using $(iii)$ of Assumption \ref{ass:fl-n}, we get
$\lim_{n \rightarrow +\infty} H_{3}^{[n]}(\bF_{n}) = \sigma^{2}.$ This ends the proof of \eqref{eq:cvV3n-fln} since $(H_{3}^{[n]}(\bF_{n}))$ is a deterministic sequence. 
\subsubsection*{Proof of \eqref{eq:cvV4n-fln}}

We set
\begin{equation*}
h_{\ell,k}^{(n)} = 2^{k - \ell} \,\cq^{p-1-(\ell+k)}\big(\cp\big(\cq^k \tilde f_{\ell,n} \otimes^2\big)\big) \quad \text{and} \quad H_{4,n} = \sum_{\ell = 0}^{p-1} \sum_{k = 0}^{p-\ell-1} h_{\ell,k}^{(n)}
\end{equation*}
in such a that $V_{4}(n) = |\GG_{n-p}|^{-1} M_{\GG_{n-p}}(H_{4,n})$. Using, \eqref{eq:geom-erg} and $(i)$ and $(ii)$ of Assumption \ref{ass:fl-n}, we get
\begin{equation*}
|h_{\ell,k}^{(n)}| \leq 2^{k-\ell} \Pp(|\Qq^{k} \tilde{f}_{\ell,n}| \otimes^{2}) \leq C 2^{k-\ell} h_{n}^{d} \alpha^{2k}.
\end{equation*}
This implies that $|H_{4,n}| \leq c_{n}$ and then that $|V_{4}(n)| \leq c_{n}$, where the sequence $(c_{n}, n \in \NN)$ is defined by
\begin{equation*}
c_{n} = C h_{n}^{d} \ind_{\{2\alpha^{2} \leq 1\}} + C h_{n}^{d} (2\alpha^{2})^{p} \ind_{\{2\alpha^{2} > 1\}}
\end{equation*}
Using \eqref{eq:fln-aS2} and the fact that $(h_{n}, n \in \NN)$ converges to 0, we get that the sequence $(c_{n}, n \in \NN)$ converges to $0$. Thus, we have that $V_{4}(n)$ is bounded by a deterministic sequence which converges to $0$. Then \eqref{eq:cvV4n-fln} follows using Remark \ref{rem3}.
\end{proof}

\begin{lem}\label{lem:cvV2nfln}
Under the Assumptions of Theorem \ref{thm:mdp-fln}, we have $V_{2}(n) \superexp 0.$
\end{lem}
\begin{proof}
Using \reff{eq:Q2-bis}, we get:
\[
V_2(n) = V_5(n)+ V_6(n),
\]
with
\begin{align*}
V_5(n) &=  |\G_n|^{-1} \sum_{i\in \G_{n-p}} \sum_{0\leq \ell<k\leq  p } 2^{p-\ell} \cq^{p-k} \big( \tilde f_{k,n} \cq^{k-\ell} \tilde f_{\ell,n}\big)(X_i),\\
V_6(n) &=     |\G_n|^{-1} \sum_{i\in \G_{n-p}} \sum_{0\leq \ell<k<  p }\sum_{r=0}^{p-k-1}  2^{p-\ell+r} \, \cq^{p-1-(r+k)}\big(\cp\big(\cq^r \tilde f_{k,n} \sot \cq  ^{k-\ell+r} \tilde f_{\ell,n} \big)\big)(X_i).
\end{align*}
First, we set
\begin{equation*}
h_{k,\ell,r}^{(n)} = 2^{r - \ell} \, \cq^{p-1-(r+k)}\big(\cp\big(\cq^r \tilde f_{k,n} \sot \cq  ^{k-\ell+r} \tilde f_{\ell,n} \big)\big) \quad \text{and} \quad H_{6,n} = \sum_{0\leq \ell<k<  p }\sum_{r=0}^{p-k-1} h_{k,\ell,r}^{(n)}
\end{equation*}
in such a way that $V_{6}(n) = |\GG_{n-p}|^{-1} M_{\GG_{n-p}}(H_{6,n}).$ Using, \eqref{eq:geom-erg} and $(i)$ and $(ii)$ of Assumption \ref{ass:fl-n}, we get
\begin{equation*}
|h_{k,\ell,r}^{(n)}| \leq C h_{n}^{d} (2\alpha^{2})^{r} \alpha^{k-\ell}.
\end{equation*}
This implies that $|H_{6,n}| \leq c_{n}$ and then that $V_{6}(n) \leq c_{n}$, where the sequence $(c_{n}, n \in \NN)$ is defined by 
$$
c_{n} = C h_{n}^{d} \ind_{\{2 \alpha^{2} \leq 1\}} + C (2\alpha^{2})^{p} h_{n}^{d} \ind_{\{2\alpha^{2} > 1\}}.
$$
Since the sequence $(c_{n}, n \in \NN)$ is deterministic and converges to 0, it follows, using Remark \ref{rem3}, that $$V_{6}(n) \superexp 0.$$
Next, for the term $V_{5}(n)$, we have for all $k > \ell$:
\begin{align*}
\big|2^{-\ell} \cq^{p-k} \big( \tilde f_{k,n} \cq^{k-\ell} \tilde f_{\ell,n}\big)\big| &\leq 2^{-\ell} \Qq^{p-k}\big(|\tilde{f}_{k,n}||\Qq^{k-\ell}\tilde{f}_{\ell,n}|\big) \leq C 2^{-\ell} h_{n}^{d/2} \alpha^{k-\ell} \Qq^{p-k}(|\tilde{f}_{k,n}|) \\ & \leq C \alpha^{p} (2\alpha)^{-\ell} \ind_{\{k = p\}} + C h_{n}^{d} (2\alpha)^{-\ell} \alpha^{k} \ind_{\{k \leq p-1\}},
\end{align*}
where we used \eqref{eq:geom-erg} for the second inequality and $(i)$ and $(ii)$ of Assumption \ref{ass:fl-n} for the second and the last inequality. Using the latter inequality in $V_{5}(n)$, we get
\begin{equation*}
|V_{5}(n)| \leq C \big(2^{-p}\ind_{\{2\alpha < 1\}} + \alpha^{p} \ind_{\{2\alpha \geq 1\}} + h_{n}^{d}\big). 
\end{equation*}
We thus have that $V_{5}(n)$ is bounded by a deterministic sequence which converges to $0$. It then follows from Remark \ref{rem3} that $$V_{5}(n) \superexp 0.$$
From the foregoing, we get the result of Lemma since $V_{2}(n) = V_{5}(n) + V_{6}(n).$ 
\end{proof}
As a direct consequence of \eqref{eq:def-V} and the Lemmas \ref{lem:cvR2nfln}, \ref{lem:cvV1nfln} and \ref{lem:cvV2nfln}, we have the following result.
\begin{lem}\label{lem:cvVnfln}
Under the Assumptions of Theorem \ref{thm:mdp-fln}, we have $V(n) \superexp \sigma^{2}.$
\end{lem}

We now study the 4th-order exponential moment condition. We stress that this condition imply in particular the exponential Lindeberg condition (condition ({\bf C3}) in Proposition \ref{prop:mdp}). We have the following result.
\begin{lem}\label{lem:lyap-fln}
Under the Assumptions of Theorem \ref{thm:mdp-fln}, we have
\[
\limsup_{n \rightarrow +\infty} \frac{1}{b_{n}^{2}} \log\PP\Big(b_{n}^{2} \sum_{i \in \GG_{n-p}} \EE[\Delta_{n,i}(\bF_{n})^{4} | \Ff_{i}] > \delta\Big) = -\infty \quad \forall \delta > 0.
\]
\end{lem}
\begin{proof}
For all $i \in \GG_{n-p}$, we have
\begin{equation}\label{eq:M4-fln}
\EE\left[\Delta_{n,i}(\bF_{n})^{4}|\Ff_{i}\right] \leq 16 (p+1)^{3} 2^{-2n} \sum_{\ell = 0}^{p} \EE_{X_{i}}\left[M_{\GG_{p-\ell}}(\tilde{f}_{\ell,n})^{4}\right],
\end{equation}
where we have  used the definition of $\Delta_{n,i}(\bF_{n})$, the inequality $(\sum_{k=0}^r a_k)^4 \leq  (r+1)^3 \sum_{k=0}^r a_k^4$ and the branching Markov property. Using \eqref{eq:M4-fln}, we get
\begin{equation}\label{eq:IsumDelta}
b_{n}^{2} \sum_{i \in \GG_{n-p}} \EE[\Delta_{n,i}(\bF_{n})^{4} | \Ff_{i}] \leq C b_{n}^{2} p^{3} 2^{-2n} \sum_{\ell = 0}^{p} \sum_{i \in \GG_{n-p}} h_{n,\ell}(X_{i}),
\end{equation}
where $h_{n,\ell}(x) = \EE_{x}[M_{\GG_{p-\ell}}(\tilde{f}_{\ell,n})^{4}].$ 
We will now prove that the right hand side of \eqref{eq:IsumDelta} converges superexponentially to $0$ at the speed $b_{n}^{2}$, that is
\begin{equation*}
\limsup_{n \rightarrow + \infty} \frac{1}{b_{n}^{2}} \log \PP\Big( C b_{n}^{2} p^{3} 2^{-2n} \Big| \sum_{\ell = 0}^{p} \sum_{i \in \GG_{n-p}} h_{n,\ell}(X_{i}) \Big| > \delta \Big) = - \infty.
\end{equation*} 

For that purpose, we will treat the case $\ell = p$, $\ell = p-1$ and finally the case $\ell \in \{0, \ldots, p-2\}$. First, we treat the case $\ell = p.$ Set $g_{p,n} = \tilde{f}_{p,n}^{4}.$ We have
\begin{equation}\label{eq:cas-l=p}
b_{n}^{2} p^{3} 2^{-2n}  \sum_{i \in \GG_{n-p}} h_{n,p}(X_{i}) = b_{n}^{2} p^{3} 2^{-2n} \sum_{i\in\GG_{n-p}} \tilde{g}_{p,n}(X_{i}) + b_{n}^{2} p^{3} 2^{-2n} \,|\GG_{n-p}| \langle \mu, g_{p, n} \rangle.
\end{equation}
Since $b_{n}^{2} p^{3} 2^{-2n} |\GG_{n-p}| \langle \mu, g_{p, n} \rangle \leq p^{3} 2^{-p} \, b_{n}^{2} \, (|\GG_{n}|h_{n}^{d})^{-1} \rightarrow 0$  as $n \rightarrow 0$, it suffices to prove that the first term of the right hand side in \eqref{eq:cas-l=p} converges superexponentially to $0$ at the speed $b_{n}^{2}$, that is, for all $\delta > 0,$
\begin{equation}\label{eq:cv-exp-l=p}
\limsup_{n \rightarrow +\infty} \frac{1}{b_{n}^{2}} \log \PP\Big( b_{n}^{2} p^{3} 2^{-2n} \, |\sum_{i\in\GG_{n-p}} \tilde{g}_{p,n}(X_{i})| > \delta \Big) = - \infty. 
\end{equation}
As in the proof of Lemma \ref{lem:cvexp-R1n}, we can prove that
\begin{equation*}
\PP\Big( b_{n}^{2} p^{3} 2^{-2n} \, |\sum_{i\in\GG_{n-p}} \tilde{g}_{p,n}(X_{i})| > \delta \Big) \leq C \exp\Big( - \frac{\delta^{2} |\GG_{n}|^{2} h_{n}^{2d}}{C p^{3} b_{n}^{2} (\delta + C p^{3} b_{n}^{2} (|\GG_{n+p}|h_{n}^{d})^{-1})}\Big).
\end{equation*}
Taking the $\log$ and dividing by $b_{n}^{2}$, we get \eqref{eq:cv-exp-l=p}.

Next, for $\ell \in \{0, \ldots, p-1\}$, we plan to prove that the quantity $b_{n}^{2} p^{3} 2^{-2n} \sum_{\ell = 0}^{p-1} \sum_{i \in \GG_{n-p}}  h_{n,\ell}(X_{i})$ is bounded by a deterministic sequence which converges to $0$. First, for $\ell = p-1$, using the branching Markov property, $(i)$ and $(ii)$ of Assumption \ref{ass:fl-n}, we have, for all $i \in \GG_{n-p,}$
\begin{equation*} 
h_{n,p-1}(X_{i}) = \EE_{X_{i}}[M_{\GG_{1}}(\tilde{f}_{p-1,n})^{4}] \leq C \Qq(\tilde{f}_{p-1,n}^{4}) \leq C h_{n}^{-d}.
\end{equation*} 
Using \eqref{eq:speed-mdp}, this implies that
\begin{equation*}
b_{n}^{2} \, p^{3} \, 2^{-2n}\sum_{i \in \GG_{n-p}} h_{n, p-1}(X_{i}) \leq C \, b_{n}^{2} \, 2^{-p} \, p^{3} (|\GG_{n}| h_{n}^{d})^{-1} \rightarrow 0 \quad \text{as $n \rightarrow +\infty.$}
\end{equation*}
%

%
Now we consider the case $\ell \in \{0, \ldots, p-2\}$. From Lemma \ref{lem:M4} with $f$ replaced by $\tilde{f}_{\ell,n}$ and $\nu$ by the Dirac mass at $X_{i}$ ($\delta_{X_{i}}$), we have
\begin{equation}\label{eq:i-lya-2}
b_{n}^{2} \, p^{3} \, 2^{-2n} \sum_{\ell = 0}^{p-2} \sum_{i \in \GG_{n-p}}  h_{n,\ell}(X_{i}) \leq b_{n}^{2} \, |\GG_{n}|^{-2} \, p^{3} \, \sum_{\ell = 0}^{p-2} \sum_{i \in \GG_{n-p}}  \sum_{j = 1}^{9} |\psi_{j,p-\ell}|(X_{i}).
\end{equation}
For all $j \in \{1, \ldots, 9\}$, we will upper bound each term of the right hand side in \eqref{eq:i-lya-2} by a deterministic sequence which converges to $0$.
\subsubsection*{\bf Upper bound of $b_{n}^{2} |\GG_{n}|^{-2} p^{3} \sum_{\ell = 0}^{p} \sum_{i \in \GG_{n-p}} |\psi_{1,p-\ell}|(X_{i})$}
Using $(i)$ of Assumption \ref{ass:fl-n}, we have
\begin{equation*}
|\psi_{1,p-\ell}| \leq C 2^{p-\ell} \Qq^{p-\ell}(f_{\ell,n}^{4}) \leq C 2^{p-\ell} h_{n}^{-d}.
\end{equation*}
Using \eqref{eq:speed-mdp}, this implies that
\begin{equation*}
b_{n}^{2} |\GG_{n}|^{-2} p^{3} \sum_{\ell = 0}^{p} \sum_{i \in \GG_{n-p}}  |\psi_{1,p-\ell}|(X_{i}) \leq C b_{n}^{2} p^{3} (|\GG_{n}| h_{n}^{d})^{-1} \rightarrow 0 \quad \text{as $n \rightarrow + \infty.$}
\end{equation*}

\subsubsection*{\bf Upper bound of $b_{n}^{2} |\GG_{n}|^{-3} p^{3} \sum_{\ell = 0}^{p} \sum_{i \in \GG_{n-p}} |\psi_{2,p-\ell}|(X_{i})$}
Using Assumption \ref{hyp:erg-unif} and $(i)$ and $(ii)$ of Assumption \ref{ass:fl-n} for the second inequality, we get
\begin{align*}
|\psi_{2,p-\ell}| &\leq C 2^{2(p-\ell)} \sum_{k = 0}^{p-\ell-1} 2^{-k} \Qq^{k}\Pp(|\Qq^{p-k-1-\ell}(\tilde{f}_{\ell,n}^{3})|\sot|\Qq^{p-\ell-k-2} (\Qq\tilde{f}_{\ell,n})|) \\
& \leq C 2^{2(p-\ell)} \sum_{k=0}^{p-\ell-1} 2^{-k} \alpha^{p-\ell-k} \, 
 \leq \, C 2^{p-\ell}\big(\ind_{\{2\alpha < 1\}} + (p-\ell)\ind_{\{2\alpha = 1\}} + (2\alpha)^{p-\ell}\ind_{\{2\alpha > 1\}}\big).
\end{align*}
Using \eqref{eq:speed-mdp} and \eqref{eq:fln-aS2}, this implies that
\begin{equation*}
b_{n}^{2} |\GG_{n}|^{-2} p^{3} \sum_{\ell = 0}^{p} \sum_{i \in \GG_{n-p}} |\psi_{2,p-\ell}|(X_{i}) \leq C  b_{n}^{2} |\GG_{n}|^{-1}\big(p^{4} \ind_{\{2\alpha \leq 1\}} + (2\alpha)^{p}\ind_{\{2\alpha > 1\}}\big) \rightarrow 0 \quad \text{as $n \rightarrow + \infty.$}
\end{equation*}

\subsubsection*{\bf Upper bound of $b_{n}^{2} |\GG_{n}|^{-2} p^{3} \sum_{\ell = 0}^{p} \sum_{i \in \GG_{n-p}} |\psi_{3,p-\ell}|(X_{i})$}
Using $(i)$ and $(ii)$ of Assumption \ref{ass:fl-n} for the second inequality, we get
\begin{equation*}
|\psi_{3,p-\ell}| \leq 2^{2(p-\ell)} \sum_{k = 0}^{p-\ell-1} 2^{-k} \Qq^{k}\Pp(\Qq^{p-\ell-k-1}(\tilde{f}_{\ell,n}^{2}) \otimes^{2}) \leq C 2^{2(p-\ell)} \sum_{k = 0}^{p-\ell-1} 2^{-k} \leq C 2^{2(p-\ell)}.
\end{equation*}
Using \eqref{eq:speed-mdp}, this implies that
\begin{equation*}
b_{n}^{2} |\GG_{n}|^{-2} p^{3} \sum_{\ell = 0}^{p} \sum_{i \in \GG_{n-p}} |\psi_{3,p-\ell}|(X_{i}) \leq  C \, b_{n}^{2} \, p^{3} \, 2^{-n + p} \rightarrow 0 \quad \text{as $n \rightarrow + \infty.$}
\end{equation*}

\subsubsection*{\bf Upper bound of $b_{n}^{2} |\GG_{n}|^{-2} p^{3} \sum_{\ell = 0}^{p} \sum_{i \in \GG_{n-p}} |\psi_{4,p-\ell}|(X_{i})$}
Using Assumption \eqref{hyp:erg-unif} and $(i)$ and $(ii)$ of Assumption \ref{ass:fl-n} for the second inequality, we get
\begin{equation*}
|\psi_{4,p-\ell}| \leq C 2^{4(p-\ell)} \Pp\big(\Pp\big(|\Qq^{p-\ell-2}\tilde{f}_{\ell,n}| \otimes^{2}\big) \otimes^{2}\big) \leq C 2^{4(p-\ell)} \alpha^{4(p-\ell-2)} h_{n}^{2d}.
\end{equation*}
Using \eqref{eq:speed-mdp} and \eqref{eq:fln-aS2}, this implies that
\begin{multline*}
b_{n}^{2} |\GG_{n}|^{-2} p^{3} \sum_{\ell = 0}^{p} \sum_{i \in \GG_{n-p}} |\psi_{4,p-\ell}|(X_{i}) \\ 
\leq  C \, (b_{n}^{2} \, p^{4} \, 2^{- n - p} h_{n}^{2d}\ind_{\{2\alpha^{2} \leq 1\}} + b_{n}^{2} p^{3} 2^{-n+p} (2\alpha^{2})^{2p} h_{n}^{2d} \ind_{\{2\alpha^{2} > 1\}}) \rightarrow 0 \quad \text{as $n \rightarrow + \infty.$}
\end{multline*}

\subsubsection*{\bf Upper bound of $b_{n}^{2} |\GG_{n}|^{-3} p^{3} \sum_{\ell = 0}^{p} \sum_{i \in \GG_{n-p}} |\psi_{5,p-\ell}|(X_{i})$}
Using Assumption \eqref{hyp:erg-unif} and $(i)$ and $(ii)$ of Assumption \ref{ass:fl-n} for the second inequality, we get
\begin{align*}
|\psi_{5,p-\ell}| &\leq C \, 2^{4(p-\ell)} \sum_{k=2}^{p-\ell-1} \sum_{r=0}^{k-1} 2^{-2k-r} \Qq^{r}\Pp\big(\Qq^{k-r-1}\big(\Pp\big(|\Qq^{p-\ell-k-1}\tilde{f}_{\ell,n}| \otimes^{2}\big)\big) \otimes^{2}\big) \\
& \leq C \, 2^{4(p-\ell)} \sum_{k = 2}^{p-\ell-1} \sum_{r=0}^{k-1} 2^{-2k-r} h_{n}^{2d} \alpha^{4(p-\ell-k)} \\ 
& \leq C \, h_{n}^{2d} \, 2^{2(p-\ell)} \big( \ind_{\{2\alpha^{2} < 1\}} + (p-\ell) \ind_{\{2\alpha^{2} =1\}} + (2\alpha^{2})^{2(p-\ell)} \ind_{\{2\alpha^{2} > 1\}}\big).
\end{align*}
Using \eqref{eq:speed-mdp} and \eqref{eq:fln-aS2}, this implies that
\begin{multline*}
b_{n}^{2} |\GG_{n}|^{-2} p^{3} \sum_{\ell = 0}^{p} \sum_{i \in \GG_{n-p}} |\psi_{5,p-\ell}|(X_{i}) \\ 
\leq  C \, (b_{n}^{2} \, p^{4} \, 2^{- n + p} h_{n}^{2d}\ind_{\{2\alpha^{2} \leq 1\}} + b_{n}^{2} p^{3} 2^{- n + p} (2\alpha^{2})^{2p} h_{n}^{2d} \ind_{\{2\alpha^{2} > 1\}}) \rightarrow 0 \quad \text{as $n \rightarrow + \infty.$}
\end{multline*}

\subsubsection*{\bf Upper bound of $b_{n}^{2} |\GG_{n}|^{-2} p^{3} \sum_{\ell = 0}^{p} \sum_{i \in \GG_{n-p}} |\psi_{6,p-\ell}|(X_{i})$}
Using Assumption \eqref{hyp:erg-unif} and $(i)$ and $(ii)$ of Assumption \ref{ass:fl-n} for the second inequality, we get
\begin{align*}
|\psi_{6,n}| &\leq C \, 2^{3(p-\ell)} \sum_{k=1}^{p-\ell-1} \sum_{r=0}^{k-1} 2^{-k-r} \Qq^{r}\Pp\big( \Qq^{k-r-1}\Pp\big(|\Qq^{p-\ell-k-1}\tilde{f}_{\ell,n}| \otimes^{2}\big) \sot \Qq^{p-\ell-r-1}(\tilde{f}^{2}_{\ell,n}) \big) \\
& \leq C \, 2^{3(p-\ell)} \sum_{k=1}^{p-\ell-1} \sum_{r=0}^{k-1} 2^{-k-r} h_{n}^{d} \alpha^{2(p-\ell-k)} \\
& \leq C \, h_{n}^{d} \, 2^{2(p-\ell)} \big(\ind_{\{2\alpha^{2} < 1\}} \, + \, (p-\ell)\, \ind_{\{2\alpha^{2} = 1\}} \, + \, (2\alpha^{2})^{p-\ell} \ind_{\{2\alpha^{2} > 1\}}\big).
\end{align*}
Using \eqref{eq:speed-mdp} and \eqref{eq:fln-aS2}, this implies that
\begin{multline*}
b_{n}^{2} |\GG_{n}|^{-2} p^{3} \sum_{\ell = 0}^{p} \sum_{i \in \GG_{n-p}} |\psi_{6,p-\ell}|(X_{i}) \\ 
\leq  C \, (b_{n}^{2} \, p^{4} \, 2^{- n + p} h_{n}^{d}\ind_{\{2\alpha^{2} \leq 1\}} + b_{n}^{2} p^{3} 2^{- n + p} (2\alpha^{2})^{2p} h_{n}^{d} \ind_{\{2\alpha^{2} > 1\}}) \rightarrow 0 \quad \text{as $n \rightarrow + \infty.$}
\end{multline*}

\subsubsection*{\bf Upper bound of $b_{n}^{2} |\GG_{n}|^{-2} p^{3} \sum_{\ell = 0}^{p} \sum_{i \in \GG_{n-p}} |\psi_{7,p-\ell}|(X_{i})$}
In the same way as for $\psi_{6,p-\ell}$, we have
\begin{multline*}
b_{n}^{2} |\GG_{n}|^{-2} p^{3} \sum_{\ell = 0}^{p} \sum_{i \in \GG_{n-p}} |\psi_{7,p-\ell}|(X_{i}) \\ 
\leq  C \, (b_{n}^{2} \, p^{4} \, 2^{- n} h_{n}^{d}\ind_{\{2\alpha \leq 1\}} + b_{n}^{2} p^{3} 2^{- n + p} (2\alpha^{2})^{2p} h_{n}^{d} \ind_{\{2\alpha > 1\}}) \rightarrow 0 \quad \text{as $n \rightarrow + \infty.$}
\end{multline*}

\subsubsection*{\bf Upper bound of $b_{n}^{2} |\GG_{n}|^{-2} p^{3} \sum_{\ell = 0}^{p} \sum_{i \in \GG_{n-p}} |\psi_{8,p-\ell}|(X_{i})$}
Using Assumption \eqref{hyp:erg-unif} and $(i)$ and $(ii)$ of Assumption \ref{ass:fl-n} for the second inequality, we get
\begin{align*}
|\psi_{8,p-\ell}| &\leq C \, 2^{4(p-\ell)} \, \sum_{k=2}^{p-\ell-1} \sum_{r=1}^{k-1} \sum_{j=0}^{r-1} 2^{-k-r-j} \Qq^{j} \Pp\big( \Qq^{r-j-1}\Pp\big(|\Qq^{p-\ell-1-r} \tilde{f}_{\ell,n}| \otimes^{2}\big) \\ 
& \hspace{5cm} \sot \Qq^{k-j-1}\Pp\big(|\Qq^{p-\ell-1-k} \tilde{f}_{\ell,n}| \otimes^{2}\big) \big) \\
&\leq C \, 2^{4(p-\ell)} \, \sum_{k=2}^{p-\ell-1} \sum_{r=1}^{k-1} \sum_{j = 0}^{r-1} 2^{-k-r-j} h_{n}^{2d} \alpha^{4(p-\ell) -2r - 2k}\\
&\leq C \, h_{n}^{2d} \, 2^{2(p-\ell)} \big(\ind_{\{2\alpha^{2} < 1\}} \, + \, (p-\ell)^{2} \, \ind_{\{2\alpha^{2} = 1\}} \, + \, (2\alpha^{2})^{2(p-\ell)} \, \ind_{\{2\alpha^{2} > 1\}}\big).
\end{align*}
Using \eqref{eq:speed-mdp} and \eqref{eq:fln-aS2}, this implies that
\begin{multline*}
b_{n}^{2} |\GG_{n}|^{-2} p^{3} \sum_{\ell = 0}^{p} \sum_{i \in \GG_{n-p}} |\psi_{8,p-\ell}|(X_{i}) \\ 
\leq  C \, (b_{n}^{2} \, p^{5} \, 2^{- n + p} h_{n}^{2d}\ind_{\{2\alpha^{2} \leq 1\}} + b_{n}^{2} p^{3} 2^{- n + p} (2\alpha^{2})^{2p} h_{n}^{2d} \ind_{\{2\alpha^{2} > 1\}}) \rightarrow 0 \quad \text{as $n \rightarrow + \infty.$}
\end{multline*}

\subsubsection*{\bf Upper bound of $b_{n}^{2} |\GG_{n}|^{-2} p^{3} \sum_{\ell = 0}^{p} \sum_{i \in \GG_{n-p}} |\psi_{9,p-\ell}|(X_{i})$}
In the same way as for $\psi_{8,p-\ell}$, we have
\begin{multline*}
b_{n}^{2} |\GG_{n}|^{-2} p^{3} \sum_{\ell = 0}^{p} \sum_{i \in \GG_{n-p}} |\psi_{9,p-\ell}|(X_{i}) \\ 
\leq  C \, (b_{n}^{2} \, p^{5} \, 2^{- n + p} h_{n}^{2d}\ind_{\{2\alpha^{2} \leq 1\}} + b_{n}^{2} p^{3} 2^{- n + p} (2\alpha^{2})^{2p} h_{n}^{2d} \ind_{\{2\alpha^{2} > 1\}}) \rightarrow 0 \quad \text{as $n \rightarrow + \infty.$}
\end{multline*}
Putting together all the upper bounds for $\ell \in \{0, \ldots, p-1\}$ and using \eqref{eq:M4-fln} and \eqref{eq:i-lya-2}, we deduce that
$b_{n}^{2} \, p^{3} \, 2^{-2n} \sum_{\ell = 0}^{p-1} \sum_{i \in \GG_{n-p}}  h_{n,\ell}(X_{i})$ is bounded by a deterministic sequence which converges to 0. As a consequence, it follows, using Remark \ref{rem3}, that
\[
\limsup_{n \rightarrow +\infty} \frac{1}{b_{n}^{2}} \log \PP\Big( b_{n}^{2} \, p^{3} \, 2^{-2n} \sum_{\ell = 0}^{p-1} \sum_{i \in \GG_{n-p}}  h_{n,\ell}(X_{i}) > \delta \Big) = - \infty.
\]
Finally, using \eqref{eq:M4-fln}, \eqref{eq:IsumDelta}, \eqref{eq:cv-exp-l=p}, we get
 \[
\limsup_{n \rightarrow +\infty} \frac{1}{b_{n}^{2}} \log\PP\Big(\sum_{i \in \GG_{n-p}} \EE[\Delta_{n,i}(\bF_{n})^{4} | \Ff_{i}] > \frac{\delta}{b_{n}^{2}}\Big) = -\infty \quad \forall \delta > 0.
\]
\end{proof}

For Chen-Ledoux type condition, we have the following result.
\begin{lem}\label{lem:C-L-fln}
Under the assumptions of Theorem \ref{thm:mdp-fln}, we have
\begin{equation*}\label{eq:chen-L}
\limsup_{n \rightarrow +\infty} \frac{1}{b_{n}^{2}} \log \Big(|\GG_{n}| \sup_{i \in \GG_{n-p}} \PP_{\Ff_{i}}\Big(|\Delta_{n,i}(\bF_{n})| > b_{n}\Big)\Big) = - \infty.
\end{equation*}
\end{lem}

\begin{proof}
For all $i \in \GG_{n-p}$, using \eqref{eq:def-DiF} we have
\begin{equation}\label{eq:ProbaCL-fln}
\PP_{\Ff_{i}}\left(|\Delta_{n,i}(\bF)| > b_{n} \sqrt{n}\right) \leq \PP_{\Ff_{i}}\left(|N_{n,i}(\bF)| > b_{n} \sqrt{n}/2 \right) + \PP_{\Ff_{i}}\left(|\EE_{X_{i}}\left[N_{n,i}(\bF)\right]| > b_{n} \sqrt{n}/2 \right),
\end{equation}
with $N_{n,i}(\bF)$ defined in \eqref{eq:def-NiFn}. Following the proof of \eqref{eq:I-R0n-CL}, we get
\begin{align*}
\PP_{\Ff_{i}}\Big(|N_{n,i}(\bF)| > \frac{b_{n} \sqrt{n}}{2} \Big) &= \PP_{X_{i}}\Big(|\sum_{\ell = 0}^{p} M_{\GG_{p-\ell}}(\tilde{f}_{\ell},n)| > \frac{b_{n}\sqrt{|\GG_{n}|}}{2}\Big) \\ 
& \leq C \exp\Big(- \frac{b_{n}^{2} |\GG_{n}|}{2(Cb_{n}|\GG_{n}|^{1/2} h_{n}^{-d/2} + 6 C |\GG_{p}|)}\Big).
\end{align*}
Next, for 
\[
\lambda = \frac{b_{n} \sqrt{|\GG_{n}|}}{2(c_{2} h_{n}^{-d/2} b_{n} \sqrt{n} + 3 c_{1} |\GG_{p}|)},
\] 
we have
\begin{align*}
\PP_{\Ff_{i}}\Big(\EE_{X_{i}}\left[N_{n,i}(\bF)\right] > \frac{b_{n} \sqrt{n}}{2} \Big) &= \PP_{X_{i}}\Big( \sum_{\ell = 0}^{p} 2^{p-\ell} \Qq^{p-\ell}(\tilde{f}_{\ell})(X_{i}) > \frac{b_{n} \sqrt{n|\GG_{n}|}}{2}\Big) \\
&\leq \exp\Big( - \frac{\lambda b_{n} \sqrt{n |\GG_{n}|}}{2} \Big) \EE_{X_{i}}\Big[ \exp\Big( \lambda \sum_{\ell = 0}^{p} 2^{p-\ell} \Qq^{p-\ell}(\tilde{f}_{\ell})(X_{i}) \Big) \Big] \\
& \leq C \exp\Big(- \frac{b_{n}^{2} |\GG_{n}|}{2(Cb_{n}|\GG_{n}|^{1/2} h_{n}^{-d/2} + 6 C |\GG_{p}|)}\Big),
\end{align*}
where we used \eqref{eq:Q1} and the branching Markov property for the first equality, Chernoff bound for the first inequality and \eqref{eq:geom-erg} for the last inequality. Doing the same thing for $-\bF$ instead of $\bF,$ we get
\begin{equation*}
\PP_{\Ff_{i}}\Big(|\EE_{X_{i}}\left[N_{n,i}(\bF)\right]| > \frac{b_{n} \sqrt{n}}{2} \Big) \leq 2 \, C \exp\Big(- \frac{b_{n}^{2} |\GG_{n}|}{2(Cb_{n}|\GG_{n}|^{1/2} h_{n}^{-d/2} + 6 C |\GG_{p}|)}\Big).
\end{equation*}
From the foregoing, we get, using \eqref{eq:ProbaCL-fln},
\begin{equation*}
|\GG_{n}| \sup_{i \in \GG_{n-p}} \PP_{\Ff_{i}}\left(|\Delta_{n,i}(\bF)| > b_{n} \sqrt{n}\right) \leq C \, |\GG_{n}| \, \exp\Big(- \frac{b_{n}^{2} |\GG_{n}|}{2(Cb_{n}|\GG_{n}|^{1/2} h_{n}^{-d/2} + 6 C |\GG_{p}|)}\Big).
\end{equation*}
Finally, taking the $\log$ and dividing by $b_{n}^{2}$ in the latter inequality, we get the result of Lemma \ref{lem:C-L-fln}.
\end{proof}

\medskip

We can  now use  Proposition \ref{prop:mdp}  to  deduce  from  Lemmas \ref{lem:cvVnfln}, \ref{lem:lyap-fln} and \ref{lem:C-L-fln} that  $\Delta_n(\bF_{n})$ satisfies a moderate deviation principle with speed $b_{n}^{2}$ and rate function $I$ defined by: $I(x) = x^{2}/(2 \sigma^{2})$ for all $x \in \RR$, with the finite variance $\sigma^{2}$ defined in \eqref{eq:limf-ln}.  Using \eqref{eq:equiexp-N-D} and Remark \ref{rem:cv-det-exp}, we then  deduce Theorem \ref{thm:mdp-fln}.

\section{Appendix}\label{sec:appendix}

%

We recall here a simplified version of Theorem 1 in \cite{djellout2002moderate}. We consider the real martingale $(M_{n},n\in\NN)$ with respect to the filtration $(\Hh_{n}, n\in \NN)$ and we denote $(\langle M\rangle_{n}, n \in \NN)$ its bracket. 

\begin{prop}\label{prop:mdp} Let $(b_{n})$ a sequence satisfying
\[b_{n} \quad \text{is increasing}, \quad b_{n} \longrightarrow +\infty, \quad \frac{b_{n}}{\sqrt{n}}\longrightarrow 0,\]
such that $c(n) := \sqrt{n}/b_{n}$ is non-decreasing, and define the
reciprocal function $c^{-1}(t)$ by
\[c^{-1}(t):=\inf\{n\in \mathbb{N}: c(n)\geq t\}.\]
Under the following conditions:
\begin{enumerate}
\item [(\textbf{C1})] there exists $Q\in \mathbb{R}_{+}^{*}$ such that for all $\delta > 0$,

$\displaystyle \limsup_{n \rightarrow +\infty} \frac{1}{b_{n}^{2}} \log\left(\PP\left(\left|\frac{\langle M\rangle_{n}}{n} - Q\right| > \delta\right)\right) = -\infty,$

\item [(\textbf{C2})]  $\displaystyle \limsup\limits_{n\rightarrow
+\infty}\frac{1}{b_{n}^{2}}\log\left(n \quad \underset{1\leq k\leq
c^{-1}(b_{n+1})}{\rm ess\,sup}
\mathbb{P}(|M_{k}-M_{k-1}| > b_{n}\sqrt{n} \Big|\mathcal{H}_{k-1})\right)=-\infty,$

\item [(\textbf{C3})] for all $ a>0$ and for all $\delta > 0$,

$\displaystyle \limsup_{n \rightarrow +\infty} \frac{1}{b_{n}^{2}} \log\left(\PP\left(\frac{1}{n}\sum\limits_{k=1}^{n}\mathbb{E}\left(|M_{k}-M_{k-1}|^{2} \mathbf{1}_{\{|M_{k}-M_{k-1}| \geq a \frac{\sqrt{n}}{b_{n}}\}} \Big|\mathcal{H}_{k-1} \right) > \delta\right)\right) = - \infty,$
\end{enumerate}
$\left(M_{n}/(b_{n}\sqrt{n})\right)_{n\in \N}$ satisfies the MDP on $\mathbb{R}$ with the
speed $b_{n}^{2}/n$ and rate function $\displaystyle I(x) =\frac{x^{2}}{2Q}.$
\end{prop}

We have the following many-to-one formulas. Ideas of the proofs can be found in \cite{Guyon} and \cite{BDG14}.
\begin{lem}\label{lem:Qi}
Let $f,g\in \cb(S)$, $x\in S$ and $n\geq m\geq 0$. Assuming that all the quantities below are well defined,  we have:
\begin{align}
\label{eq:Q1} \E_x\left[M_{\G_n}(f)\right] &=|\G_n|\, \cq^n f(x)= 2^n\, \cq^n f(x) ,\\
\label{eq:Q2} \E_x\left[M_{\G_n}(f)^2\right] &=2^n\, \cq^n (f^2) (x) + \sum_{k=0}^{n-1} 2^{n+k}\,   \cq^{n-k-1}\left( \cp \left(\cq^{k}f\otimes \cq^k f \right)\right) (x),\\
\label{eq:Q2-bis} \E_x\left[M_{\G_n}(f)M_{\G_m}(g)\right] &=2^{n} \cq^{m} \left(g \cq^{n-m} f\right)(x)\\
\nonumber &\hspace{2cm} + \sum_{k=0}^{m-1} 2^{n+k}\, \cq^{m-k-1} \left(\cp\left(\cq^k g \sot \cq^{n-m+k} f\right) \right)(x). 
\end{align}
\end{lem}

%
We recall  the following result due to Bochner (see \cite[Theorem 1A]{Parzen1962}  which can be easily extended to any dimension $d\geq 1$). 
\begin{lem}\label{lem:bochner}
Let $(h_{n},n\in\NN)$ be a sequence of positive numbers converging to $0$ as $n$ goes to infinity. Let $g: \RR^{d} \rightarrow \RR$ be a measurable function such that $\int_{\RR^{d}} |g(x)|dx < +\infty$. Let $f: \RR^{d} \rightarrow \RR$  be a measurable function such that  $\norm{f}_{\infty}<+\infty $, $\int_{\R^d} |f(y)|\, dy < + \infty$  and $\lim_{|x|\rightarrow +\infty} |x|f(x)=0$. Define
\begin{equation*}
g_{n}(x) = h_{n}^{-d}\int_{\RR^{d}} f(h_{n}^{-1}(x-y))g(y)dy.
\end{equation*}
Then, we have at every point $x$ of continuity of $g$,
\begin{equation*}
\lim_{n\rightarrow +\infty} g_{n}(x) = g(x)\int_{\RR} f(y)dy.
\end{equation*} 
\end{lem}
We also give some bounds on $\E_x\left[M_{\G_{n}}(f) ^4\right]$,  see the proof of Theorem 2.1 in \cite{BDG14}.  We will use the notation:
\[
g\otimes^2=g\otimes g.
\]
\begin{lem}\label{lem:M4}
There exists a  finite constant $C$ such that for  all $f\in \cb(S)$,$n\in \N$ and  $\nu$ a probability measure on $S$,  assuming that all the quantities below are well defined, there exist functions $\psi_{j,n}$ for $1\leq j\leq 9$ such that:
\[
\E_\nu\left[M_{\G_n}(f)^4\right]= \sum_{j=1}^9 \langle \nu, \psi_{j, n}\rangle,
\]
and, with $h_{k}= \cq^{k - 1} (f) $ and (notice that either $|\psi_j|$ or $|\langle \nu, \psi_j \rangle|$ is bounded), writing  $\nu   g = \langle \nu  ,   g   \rangle$:
\begin{align*}
| \psi_{1, n}| &\leq C \,2^n \cq^n(f^4),\\
| \nu \psi_{2, n}| &\leq C\,  2^{2n}\, \sum_{k=0}^{n-1} 2^{-k} |\nu \cq^k \cp  \left( \cq^{n-k - 1}( f^3) \sot h_{n- k} \right)|,\\
|\psi_{3, n}| &\leq C 2^{2n} \sum_{k=0}^{n-1} 2^{-k}\,  \cq^k \cp \left( \cq^{n-k - 1} (f^2) \otimes^2\right),\\
|\psi_{4, n}| &\leq C \, 2^{4n} \, \cp \left( |\cp(h_{n-1}\otimes^2)\otimes^2|\right), \\
|\psi_{5, n}| &\leq C\,  2^{4n} \, \sum_{k=2}^{n-1} \sum_{r=0}^{k -1}  2^{-2k-r }  \cq^r \cp \left( \cq^{k -r- 1} |\cp (h_{n- k} \otimes^2)|\otimes^2 \right),\\
|\psi_{6, n}| &\leq C\, 2^{3n} \, \sum_{k=1}^{n-1} \sum_{r=0}^{k -1} 2^{-k-r }  \cq^r| \cp \left(\cq^{k -r-1}\cp \left( h_{n-k} \otimes^2 \right)\sot\cq^{n-r-1}(f^2) \right)|,\\
|\nu \psi_{7, n}| &\leq  C\, 2^{3n} \,   \sum_{k=1}^{n-1} \sum_{r=0}^{k -1} 2^{-k-r } |\nu \cq^r \cp \left(\cq^{k -r-1}\cp \left( h_{n-k} \sot  \cq^{n-k -1} (f^2) \right)\sot h_{n-r} \right)|,\\
|\psi_{8, n}| &\leq C\,  2^{4n} \, \sum_{k=2}^{n-1} \sum_{r=1}^{k -1} \sum_{j=0}^{r-1} 2^{-k-r-j } \cq^j \cp \left(|\cq^{r-j-1}\cp \left( h_{n-r} \otimes^2 \right)|\sot |\cq^{k-j-1}\cp \left( h_{n-k} \otimes^2 \right)|\right),\\
|\psi_{9, n}| & \leq C\,  2^{4n} \, \sum_{k=2}^{n-1} \sum_{r=1}^{k -1} \sum_{j=0}^{r-1}2^{-k-r-j } \cq^j |\cp \left(\cq^{r-j-1}|\cp \left( h_{n-r} \sot \cq^{k -r -1}\cp\left(h_{n-k}\otimes^2 \right)\right) \sot h_{n-j} \right)|.
\end{align*}
\end{lem}

\bibliographystyle{abbrv}
\bibliography{biblio}
\end{document}